\documentclass[11pt,leqno]{article}

\usepackage{mathrsfs}
\usepackage{amsmath}
\usepackage{amssymb}
\usepackage{amsthm}

\numberwithin{equation}{section}

\title{On the Gevrey strong hyperbolicity}

\author{ Tatsuo Nishitani}

 \date{}
\begin{document}

\maketitle

\begin{abstract}{In this paper we are concerned with a homogeneous differential operator $p$ of order $m$ of which characteristic set of order $m$ is assumed to be a smooth manifold. We define the Gevrey strong hyperbolicity index as the largest number $s$ such that the Cauchy problem for $p+Q$ is well-posed in the Gevrey class of order $s$ for any differential operator $Q$ of order less than $m$. We study the case of the largest index and we discuss in which way the Gevrey strong hyperbolicity index relates with the geometry  of bicharacteristics of $p$ near the characteristic manifold.}
\end{abstract}


\maketitle

\def\R{{\mathbb R}}
\def\C{{\mathbb C}}
\def\lr#1{\langle{#1}\rangle}
\def\al{\alpha}
\def\be{\beta}
\def\ga{\gamma}
\def\si{\sigma}
\def\la{\lambda}
\def\La{\Lambda}
\def\dif{\partial}
\def\N{{\mathbb N}}
\def\xim{\langle{\xi}\rangle_{\gamma}}
\def\xig{\langle{\xi'}\rangle_{\gamma}}
\def\bg{{\bar g}}

\section{Introduction}

Let 
\[
P=D_0^m+\sum_{|\al|\leq m, \al_0<m}a_{\al}(x)D^{\al}=p(x,D)+P_{m-1}(x,D)+\cdots
\]
be a differential operator of order $m$ defined near the origin of $\R^{n+1}$ where $x=(x_0,\ldots,x_n)=(x_0,x')$ and 
\[
D_j=-i\dif/\dif x_j,\;\;D=(D_0,D'),\;\;D'=(D_1,\ldots,D_n).
\]
Here $p(x,\xi)$ is the principal symbol of $P$;
\[
p(x,\xi)=\xi_0^m+\sum_{|\alpha|=m,\al_0<m}a_{\al}(x)\xi^{\al}.
\]
We assume that the coefficients $a_{\al}(x)$ are in the Gevrey class of order $s>1$, sufficiently close to $1$, which are constant outside $|x'|\leq R$.  We say that $f(x)\in \gamma^{(s)}(\R^{n+1})$, the Gevrey class of order $s$, if for any compact set $K\subset \R^{n+1}$ there exist $C>0, A>0$ such that we have
\[
|D^{\al}f(x)|\leq CA^{|\al|}|\al|!^s,\;\;x\in K,\;\;\forall \alpha\in\N^{n+1}.
\]
\newtheorem{definition}{ Definition}[section]
\begin{definition}\rm
\label{def:solvability}
We say that the Cauchy problem for $P$ is $\gamma^{(s)}$ well-posed at the origin if  for any $\Phi=(u_0,u_1,\ldots,u_{m-1})\in (\gamma^{(s)}(\R^n))^m$ there exists a neighborhood $U_{\Phi}$ of the origin such that the Cauchy problem 
\[
\left\{\begin{array}{l}
Pu=0\quad \mbox{in}\quad U_{\Phi},\\
D_0^ju(0,x')=u_j(x'),\; j=0,1,\ldots,m-1,\;\; x'\in U_{\Phi}\cap\{x_0=0\}
\end{array}\right.
\]
has a unique solution $u(x)\in C^{\infty}(U_{\Phi})$. 
\end{definition}
It is a fundamental fact  that if $p(x,\xi)$ is  strictly hyperbolic near the origin, that is $p(x,\xi_0,\xi')=0$ has $m$ real distinct roots for any $x$, near the origin and any $\xi'\neq 0$ then  the Cauchy problem for  $p+Q$ with any differential operator $Q$ of order less than $m$ is $C^{\infty}$ well-posed near the origin. In particular, $\gamma^{(s)}$ well-posed for any $s>1$. On the other hand the Lax-Mizohata theorem in the Gevrey classes asserts:
\newtheorem{pro}{Proposition}[section]
\begin{pro}[{\cite[Theorem 2.2]{Ni0}}] If the Cauchy problem for $P$ is $\gamma^{(s)}$~{\rm ($s>1$)} well-posed at the origin then $p(0,\xi_0,\xi')=0$ has only real roots $\xi_0$ for any $\xi'\in\R^n$.
\end{pro}
Taking this result into account we assume, throughout the paper,  that $p(x,\xi_0,\xi')=0$ has only real roots for any $x$ near the origin and any $\xi'\in \R^n$.

\begin{definition}
\label{def:index} \rm We define $G(p)$ {\rm (}the Gevrey strong hyperbolicity index {\rm )} by
\[
G(p)=\sup\Big\{1\leq s\Big| \begin{array}{ll}\text{Cauchy problem for $p+Q$ is $\gamma^{(s)}$ well-posed at the }\\
\text{origin for any differential operator $Q$ of order $<m$}
\end{array}\Big\}.
\]
\end{definition}

We first recall a basic result of Bronshtein \cite{Br}.
\newtheorem{theorem}{Theorem}[section]
\begin{theorem}[{\cite[Theorem 1]{Br}}]
\label{thm:Bron} Let $p$ be a homogeneous differential operator of order $m$ with real characteristic roots. Then for any differential operator $Q$ of order less than $m$, the Cauchy problem for $p+Q$ is $\gamma^{(m/(m-1))}$ well-posed.
\end{theorem}
This implies that for differential operators $p$ of order $m$ with real characteristic roots   we have
\[
G(p)\geq m/(m-1).
\]
We also recall a result which bounds $G(p)$ from above. The following result is a special case of Ivrii \cite[Theorem 1]{Iv0}. Recall that $(x,\xi)\in \R^{n+1}\times(\R^{n+1}\setminus\{0\})$ is called a characteristic of order $r$ of $p$ if
\[
\dif_x^{\al}\dif_{\xi}^{\be}p(x,\xi)=0, \quad \forall |\al+\be|<r.
\]
\begin{theorem}[{\cite[Theorem 1]{Iv0}}]
\label{thm:Ivrii} Let $p$ be a homogeneous differential operator of order $m$ with real analytic coefficients and let $(0,{\bar\xi})$, ${\bar\xi}=(0,\ldots,0,1)\in \R^{n+1}$ be a characteristic of order $m$. If the Cauchy problem for $P=p+P_{m-1}+\cdots$ is $\gamma^{(\kappa)}$ well-posed at the origin we have
\[
\dif_{\xi}^{\alpha}\dif_x^{\beta}P_{m-1}(0,{\bar\xi})=0
\]
for any $|\alpha+\beta|\leq  m-2\kappa/(\kappa-1)$.
\end{theorem}
Assume that $p$ has a characteristic $(0,{\bar\xi})$ of order $m$ and that the Cauchy problem for $p+P_{m-1}+\cdots$ is  $\gamma^{(\kappa)}$ well-posed  for any $P_{m-1}$. Then  from Theorem \ref{thm:Ivrii} it follows that $m-2\kappa/(\kappa-1)< 0
$, that is $\kappa<m/(m-2)$ which yields
\[
G(p)\leq m/(m-2).
\]

Let $\rho$ be a characteristic of order $m$. Then the localization $p_{\rho}(X)$ of $p$ at $\rho$ is defined by $p(\rho+\mu X)=\mu^m(p_{\rho}(X)+o(1))$ with $X=(x,\xi)$ as $ \mu\to 0$ which is nothing but the first non-vanishing term of the Taylor expansion of $p$ around $\rho$. Note that $p_{\rho}$ is a hyperbolic polynomial in $X$ in the direction $(0,\theta) \in\R^{n+1}\times \R^{n+1}$ where $\theta=(1,\ldots,0)\in\R^{n+1}$ (for example \cite[Lemma 8.7.2]{Hobook:iti}). The hyperbolic cone $\Gamma_{\rho}$ of $p_{\rho}$ is the connected component of $(0,\theta)$ in the set (for example \cite[Lemma 8.7.3]{Hobook:iti})
\[
\Gamma_{\rho}=\{X\in \R^{2(n+1)}\mid p_{\rho}(X)\neq 0\}
\]
and the propagation cone $C_{\rho}$ of the localization $p_{\rho}$ is 
given by
$$
C_{\rho}=\{X\in \R^{2(n+1)}\mid (d\xi\wedge dx)(X,Y)\leq 0,\forall Y\in \Gamma_{\rho}\}.
$$
Let $H_p=\sum_{j=0}^n(\partial p/\partial \xi_j)\partial/\partial x_j-(\partial p/\partial x_j )\partial/\partial \xi_j$
be the Hamilton vector field of $p$. The integral curves of $H_p$, along which $p=0$, are called bicharacteristics of $p$. 
We note that $C_{\rho}$ is the {\it minimal} cone including every 
bicharacteristic which has $\rho$ as a limit point in the following sense:
\newtheorem{lemma}{Lemma}[section]
\begin{lemma}[{\cite[Lemma 1.1.1]{KN}}]
\label{lem:honeone} Let $\rho\in \R^{n+1}\times(\R^{n+1}\setminus \{0\})$ be a multiple characteristic 
of $p$. Assume that there are simple characteristics $\rho_j$ and non-zero real numbers  
$\ga_j$ with $\gamma_jp_{\rho_j}(0,\theta)>0$ such that
$$
\rho_j\to \rho\mbox{~~and~~}\ga_j H_p(\rho_j)\to X,\quad j\to \infty.
$$
Then $X\in C_{\rho}$.
\end{lemma}
We now introduce assumptions of which motivation will be discussed in the next section. Denote by $\Sigma$ the set of characteristics of order $m$ of $p(x,\xi)$;
\[
\Sigma=\{(x,\xi)\in \R^{n+1}\times (\R^{n+1}\setminus \{0\})\mid \dif_x^{\al}\dif_{\xi}^{\be}p(x,\xi)=0, \forall |\al+\be|<m\}
\]
which is assumed to be a $\gamma^{(s)}$ manifold. 
Note that $p_{\rho}$ is a function on $\R^{2(n+1)}/T_{\rho}\Sigma$ because $p_{\rho}(X+Y)=p_{\rho}(Y)$ for any $X\in T_{\rho}\Sigma$ and any $Y\in\R^{2(n+1)}$. We assume that
\begin{equation}
\label{eq:assump}
p_{\rho}\;\;\text{is a strictly hyperbolic polynomial on}\;\;\R^{2(n+1)}/T_{\rho}\Sigma,\;\;\rho\in\Sigma.
\end{equation}
We also assume that the propagation cone $C_{\rho}$ is transversal to the characteristic manifold $\Sigma$;
\begin{equation}
\label{eq:odan}
C_{\rho}\cap T_{\rho}\Sigma=\{0\},\;\;\rho\in\Sigma.
\end{equation}
Our aim in this paper is to prove
\begin{theorem}
\label{thm:new}
Assume  \eqref{eq:assump} and \eqref{eq:odan}. Then the Cauchy problem for $p+Q$ is $\gamma^{(s)}$ well-posed at the origin for any differential operator $Q$ of order less than $m$ and for any $1<s<m/(m-2)$. In particular we have $G(p)=m/(m-2)$.
\end{theorem}
\newtheorem{examp}{Example}[section]
\begin{examp}
\label{ex:tyusho}\rm Let 
\[
q(\zeta)=\zeta_0^m+\sum_{|\alpha|=m, \alpha_0\leq m-2}c_{\alpha}\zeta^{\alpha},\quad \zeta=(\zeta_0,\zeta_1,\ldots,\zeta_k)
\]
be a strictly hyperbolic polynomial in the direction $\zeta_0$ where $k\leq n$.  Let $b_j(x,\xi')$, $j=1,\ldots,k$ be smooth functions in a conic neighborhood of $(0,{\hat \xi'})$ which are homogeneous of degree $1$ in $\xi'$ with linearly independent differentials at $(0,{\hat \xi}')$. We define
\[
p(x,\xi)=q(b(x,\xi)),\;\;b=(b_0,b_1,\ldots,b_k)
\]
where we set $b_0(x,\xi)=\xi_0$ for notational convenience. Then it is easy to see that $p(x,\xi)$ verifies the condition \eqref{eq:assump} near $\rho=(0,0,{\hat \xi}')$ with $\Sigma=\{(x,\xi)\mid b_j(x,\xi)=0, j=0,\ldots,k\}$ and $p_{\rho}(x,\xi)=q(db_{\rho}(x,\xi))$, that is
\[
p_{\rho}(x,\xi)=q({\hat b}(x,\xi)),\quad {\hat b}=({\hat b}_0,{\hat b}_1,\ldots,{\hat b}_k)
\]
where ${\hat b}_j(x,\xi)$ is the linear part of $b_j(x,\xi)$ at $\rho$. Therefor $\Gamma_{\rho}=\{X\mid {\hat b}(X)\in \Gamma\}$ where $\Gamma$ is the hyperbolic cone of $q$. If 
\begin{equation}
\label{eq:onaa}
\big(\{b_i,b_j\}\big)_{0\leq i,j\leq k}\quad\mbox{is non-singular at $\rho$}
\end{equation}
 then $p(x,\xi)$ verifies the condition \eqref{eq:odan} near $\rho$ where $\{b_i,b_j\}$ denotes the Poisson bracket
\[
\sum_{\mu=0}^n(\partial b_i/\partial \xi_{\mu})(\partial b_j/\dif x_{\mu})-(\dif b_i/\dif x_{\mu})(\dif b_j/\dif \xi_{\mu}).
\]
 We check that \eqref{eq:onaa} implies \eqref{eq:odan}. Note that \eqref{eq:odan} is equivalent to $\Gamma_{\rho}\cap (T_{\rho}\Sigma)^{\sigma}\neq \emptyset$. Since $(T_{\rho}\Sigma)^{\sigma}$ is  spanned by $H_{b_0}(\rho),H_{b_1}(\rho),\ldots,H_{b_k}(\rho)$ it suffices to show that there are $c_j$ such that  $0\neq X=\sum_{j=0}^kc_jH_{b_j}(\rho)\in \Gamma_{\rho}$. From
\[
{\hat b}_j(X)=\sum_{i=0}^kc_i\{b_i,b_j\}(\rho),\;\;j=0,\ldots,k
\]
one can choose $c_j$ so that ${\hat b}(X)=(1,0,\ldots,0)$ by assumption \eqref{eq:onaa} and hence the result.
\end{examp}
\begin{examp}
\label{ex:kihon}\rm Consider
\begin{equation}
\label{eq:kihon:q}
q(\zeta_0,\zeta_1,\zeta_2)=\prod_{j=1}^{\ell}\big(\zeta_0-c_j(\zeta_1^2+\zeta_2^2)\big)
\end{equation}
where $c_j$ are real positive constants different from each other and $2\ell=m$. Take $b_1=(x_0-x_1)\xi_n$, $b_2=\xi_1$ and consider
\[
p(x,\xi)=\prod_{j=1}^{\ell}\big(\xi_0^2-c_j((x_0-x_1)^2\xi_n^2+\xi_1^2)\big)
\]
in a conic neighborhood of $\rho=(0,0,\ldots,0,1)$. The $3\times 3$ anti-symmetric matrix $(\{b_i,b_j\})$ is obviously singular. If $\max\{c_j\}=c<1$ then $C_{\rho}\cap T_{\rho}\Sigma=\{0\}$. To see this take any $X=(t,t,x_2,\ldots,x_n,0,0,\xi_2,\ldots,\xi_n)\in T_{\rho}\Sigma$. Assume $X\in C_{\rho}$ so that $(d\xi\wedge dx)(X,Y)\leq 0$ for any $Y=(y,\eta)\in \Gamma_{\rho}$, that is for any $(y,\eta)\in \R^{2(n+1)}$ with $\eta_0^2>c((y_0-y_1)^2+\eta_1^2)$ and $\eta_0>0$. This implies that $x_2=\cdots=x_n=0$, $\xi_2=\cdots=\xi_n=0$ and $-t(\eta_0+\eta_1)\leq 0$ for any $\eta_0>\sqrt{c}\,|\eta_1|$. Since $c<1$ this gives $t=0$ so that $X=0$. 

On the other hand if $\max\{c_j\}=c\geq 1$ then $C_{\rho}\cap T_{\rho}\Sigma\neq \{0\}$. Indeed let $X=(1,1,0,\ldots,0,0,\ldots,0)\in T_{\rho}\Sigma$. Noting that $\eta_0>\sqrt{c}\,|\eta_1|$ if $Y=(y,\eta)\in \Gamma_{\rho}$ we see $(d\xi\wedge dx)(X,Y)=-\eta_0-\eta_1<0$ for any $Y\in \Gamma_{\rho}$  which proves $X\in C_{\rho}$.
\end{examp}
\begin{examp}
\label{ex:kihon:2}\rm Take $q$ in \eqref{eq:kihon:q} and choose $b_1=x_0\xi_n$, $b_2=\xi_1$ and consider 
\[
p(x,\xi)=\prod_{j=1}^{\ell}\big(\xi_0^2-c_j(x_0^2\xi_n^2+\xi_1^2)\big)
\]
 near $\rho=(0,0,\ldots,0,1)$. As remarked in Example \ref{ex:kihon} the matrix $(\{b_i,b_j\})$ is singular. Suppose  $X=(0,x_1,\ldots,x_n,0,0,\xi_2,\ldots,\xi_n)\in T_{\rho}\Sigma\cap C_{\rho}$. As in Example \ref{ex:kihon} we conclude $x_2=\cdots=x_n=0$, $\xi_1=\cdots=\xi_n=0$ and $-x_1\eta_1\leq 0$ for any $\eta_0^2>c(y_0^2+\eta_1^2)$. This gives $x_1=0$ so that $X=0$. Thus we conclude $C_{\rho}\cap T_{\rho}\Sigma=\{0\}$.
\end{examp}
\begin{examp}
\label{ex:CT}\rm We specialize Example \ref{ex:tyusho} with 
\[
q(\zeta_0,\zeta_1)=\prod_{j=1}^m(\zeta_0-\al_j\zeta_1),\quad q(\zeta_0,\zeta_1)=\prod_{j=1}^{\ell}(\zeta_0^2-c_j\zeta_1^2)
\]
where $\al_j$ are real constants different from each other such that $\sum_{j=1}^m\al_j=0$ and $c_j$ are positive constant different from each other and $m=2\ell$. For these $q$ choosing $b_1=x_0\xi_1$ and  $b_1=x_0|\xi'|$ respectively we get
\[
p(x,\xi)=\prod_{j=1}^m(\xi_0-\al_j x_0\xi_1),\quad p(x,\xi)=\prod_{j=1}^{\ell}(\xi_0^2-c_jx_0^2|\xi'|^2).
\]
It is clear that $\{b_0,b_1\}=\xi_1\neq 0$ and $\{b_0,b_1\}=|\xi'|\neq 0$ respectively and hence $C_{\rho}\cap T_{\rho}\Sigma=\{0\}$.  We find these examples in  \cite{CT} where they studied Levi type conditions for differential operators of order $m$ with coefficients depending only on the time variable. 
\end{examp}
\section{Motivation, the doubly characteristic case}

In this section we provide the motivation to introduce  $G(p)$ and  assumptions  \eqref{eq:assump}, \eqref{eq:odan}. Let $m=2$ and we consider differential operators of second order
\[
P(x,D)=p(x,D)+P_1(x,D)+P_0(x)
\]
of principal symbol $p(x,\xi)$. Let $\rho$ be a double characteristic of $p$ and hence  singular (stationary) point of $H_p$. We linearize the Hamilton equation ${\dot X}=H_p(X)$ at $\rho$, the linearized equation turns to be ${\dot Y}=F_p(\rho)Y$ where $F_p(\rho)$ is given by
\[
F_p(\rho)=\left(\begin{array}{cc}
\displaystyle{\frac{\partial^2 p}{\partial x\partial \xi}(\rho)}&\displaystyle{\frac{\partial^2 p}{\partial \xi\partial \xi}(\rho)}\\[8pt]
\displaystyle{-\frac{\partial^2 p}{\partial x\partial x}(\rho)}&\displaystyle{-\frac{\partial^2 p}{\partial \xi\partial x}(\rho)}\end{array}\right)
\]
and called  the Hamilton map {\rm(}fundamental matrix{\rm)} of $p$ at $\rho$.

The special structure of $F_p(\rho)$ results from the fact that  $p(x,\xi_0,\xi')=0$ has only real roots $\xi_0$ for any $(x,\xi')$.
\begin{lemma}[{\cite[Lemma 9.2, 9.4]{IP}}]
\label{lem:1:a} All eigenvalues of the Hamilton map $F_p(\rho)$ are on the imaginary axis, possibly one exception of a pair of non-zero real eigenvalues.
\end{lemma}
We assume that the doubly characteristic set $\Sigma=\{(x,\xi)\mid \dif_{\xi}^{\al}\dif_x^{\be}p(x,\xi)=0, \forall  |\al+\be|<2\}$ verifies the following conditions:
\begin{equation}
\label{eq:fundassum}
\begin{cases}
\mbox{$\Sigma$ is a $\gamma^{(s)}$ manifold},\\
\text{$p$ vanishes on $\Sigma$ of order exactly $2$},\\
{\rm rank }\,(d\xi\wedge dx)=\text{const. on}\;\;\Sigma.
\end{cases}
\end{equation}
Note that
 $p_{\rho}(X)$ is {\it always a strictly hyperbolic polynomial } on $\R^{2(n+1)}/T_{\rho}\Sigma$ as far as $p$ vanishes on $\Sigma$ of order exactly $2$. We also assume that the codimension $\Sigma$ is $3$ and no transition of spectral type of $F_p$ occur on $\Sigma$, that is we assume
\begin{equation}
\label{eq:second}
\text{either}\quad {\rm Ker}\,F^2_p\cap {\rm Im}\,F^2_p=\{0\}\;\;\text{or}\;\;{\rm Ker}\,F^2_p\cap {\rm Im}\,F^2_p\neq\{0\}
\end{equation}
throughout $\Sigma$. The following table sums up a general picture of the Gevrey strong hyperbolicity for differential operators with double characteristics (\cite{BN4, BN3, Ni11, Iwa2}) where $W={\rm Ker}\,F^2_p\cap {\rm Im}\,F^2_p$.

\medskip
\begin{tabular}{|p{2.8cm}|p{1.6cm}|p{4cm}|p{1.8cm}|}

\hline

Spectrum of $F_p$&\vspace{1mm}$W$&Geometry of bicharacteristics near $\Sigma$&\vspace{1mm}$G(p)$\\

\hline

\hline

Exists non-zero real eigenvalue&\vspace{1mm}$W=\{0\}$&At every point on $\Sigma$ exactly two bicharacteristics intersect $\Sigma$ transversally&\vspace{2mm}$G(p)=\infty$\\

\hline

No non-zero real eigenvalue&\vspace{4mm}$W\neq \{0\}$&No bicharacteristic intersects  $\Sigma$&\vspace{1mm}$G(p)=4$\\

      \cline{3-4}

      &&Exists a bicharacteristic tangent to 
      $\Sigma$&\vspace{1mm}$G(p)=3$\\

      \cline{2-4}

      &\vspace{0.5mm}$W=\{0\}$&No bicharacteristic intersects  $\Sigma$&\vspace{1mm}$G(p)=2$\\

      \hline

      \end{tabular}

\medskip

This table shows that, assuming \eqref{eq:fundassum}, \eqref{eq:second} and the codimension $\Sigma$ is $3$, the Gevrey strong hyperbolicity index $G(p)$ takes only the values $2$, $3$, $4$ and $\infty$ and that these values completely determine the structure  of the Hamilton map and the geometry of bicharacteristics  near $\Sigma$ and vice versa.
\begin{lemma}[{\cite[Corollary 1.4.7]{Ho2}, \cite[Lemma 1.1.3]{KN}}]\label{lem:honethree} Let $\rho$ be a double characteristic. Then the following two conditions are equivalent.
\begin{description}
\item{\rm(i)} $F_p(\rho)$ has non-zero real eigenvalues,
\item{\rm(ii)} $C_{\rho}\cap T_{\rho}\Sigma=\{0\}$.
\end{description}
\end{lemma}
Note that the condition (ii) is well defined for characteristics of {\it any order} while $F_p(\rho)\equiv 0$ if $\rho$ is a characteristic of order larger than $2$.

\newtheorem{rem}{\rm Remark}[section]
\begin{rem}
\label{rem:yoso}\rm 
Based on the table, it is quite natural to ask whether the converse of Theorem \ref{thm:new} is true. That is if $G(p)=m/(m-2)$ then \eqref{eq:assump} and \eqref{eq:odan} hold?
\end{rem}
\begin{rem}
\label{rem:yoso:2}\rm
Consider the case $C_{\rho}\subset T_{\rho}\Sigma$ that would be considered as a opposite case to $C_{\rho}\cap T_{\rho}\Sigma=\{0\}$. Here we note 
\begin{lemma}[{ \cite[Lemma 2.11]{Ni:2} }]We have $C_{\rho}\subset T_{\rho}\Sigma$ if and only if $T_{\rho}\Sigma$ is involutive, that is $(T_{\rho}\Sigma)^{\sigma}\subset T_{\rho}\Sigma$ where $(T_{\rho}\Sigma)^{\sigma}=\{X\in \R^{n+1}\times\R^{n+1}\mid (d\xi\wedge dx)(X,Y)=0, \forall Y\in T_{\rho}\Sigma\}$.
\end{lemma}
It is also natural to ask whether $G(p)=m/(m-1)$ if $C_{\rho}\subset T_{\rho}\Sigma$, $\rho\in\Sigma$. When $\Sigma$ is involutive one can choose homogeneous symplectic coordinates $x$, $\xi$ in a conic neighborhood of $\rho\in\Sigma$ such that $\Sigma$ is defined by (\cite[Theorem 21.2.4]{Hobook}, for example)
 \[
 \xi_0=\xi_1=\cdots=\xi_k=0.
 \]
 Thus by conjugation of a Fourier integral operator $p(x,\xi)$ can be written
 \[
 p(x,\xi)=\xi_0^m+\sum_{\al_0\leq m-2,|\al|=m}a_{\al}(x,\xi){\tilde \xi}^{\al}
 \]
 where ${\tilde \xi}=(\xi_0,\xi_1,\ldots,\xi_k)$. Thus $\dif_{\xi}^{\al}\dif_x^{\be}p(0,{\bar\xi})=0$ for $|\al|<m$ and any $\be$. If the resulting $p(x,D)$ is a {\it differential operator} so that $a_{\al}(x,\xi)=a_{\al}(x)$ then from \cite[Theorem 1]{Iv0} we conclude that if the Cauchy problem for $p+P_{m-1}+\cdots$ is $\gamma^{(\kappa)}$ well-posed then $\dif_{\xi}^{\al}\dif_x^{\be}P_{m-1}(0,{\bar\xi})=0$ for any $|\al|\leq m-\kappa/(\kappa-1)$ and any $\be$. This proves $G(p)\leq m/(m-1)$ and hence $G(p)=m/(m-1)$. 
\end{rem}
\begin{examp}
\label{rem:takusan} \rm When $m\geq 3$ the geometry of $p$ with the limit point $\rho$ becomes to be complicated comparing with the case $m=2$, even \eqref{eq:assump} and \eqref{eq:odan} are satisfied. We give an example. Let us consider
\[
p(x,\xi)=\xi_0^3-3a\{(x_0^2+x_1^2)\xi_n^2+\xi_1^2\}\xi_0-2bx_0x_1\xi_1\xi^2_n
\]
near $\rho=(0,\ldots,0,1)$ which is obtained from Example \ref{ex:kihon} with
\[
q(\zeta_0,\zeta_1,\zeta_2,\zeta_3)=\zeta_0^3-3a(\zeta_1^2+\zeta_2^2+\zeta_3^2)\zeta_0-2b\zeta_1\zeta_2\zeta_3
\]
and $b_1=x_0\xi_n$, $b_2=x_1\xi_n$, $b_3=\xi_1$  where $a>0$, $b$ are real constants. Choosing $b=\delta a^{3/2}$ with $|\delta|<1$ and repeating similar arguments as in Example \ref{ex:kihon:2} it is easily seen that $p(x,\xi)$ satisfies \eqref{eq:assump} and \eqref{eq:odan}. 

Consider the Hamilton equations
\begin{equation}
\label{eq:coco}
{\dot x}_j=\partial p/\partial \xi_j,\quad {\dot \xi}_j=- \partial p/\partial x_j,\quad j=0,\ldots,n.
\end{equation}
 Since ${\dot \xi}_n=0$ we take $\xi_n=1$ and $x_1=\xi_0=0$ in \eqref{eq:coco} so that the resulting equations reduce to:
\begin{equation}
\label{eq:tezan}
{\dot x}_0=-3a(x_0^2+\xi_1^2),\quad {\dot \xi}_1=2bx_0\xi_1.
\end{equation}
We fix $-1<\delta<0$ and take $a>0$ so that $2b/(3a)<-1$. Then any integral curve of \eqref{eq:tezan} passing a point in the cone $|\xi_1|<|1+(2b/3a)|^{1/2}|x_0|$, $x_0<0$ arrives at the origin inside the cone (see, for example \cite{Tez}). In particular there are infinitely many bicharacteristics with the limit point $\rho$.
\end{examp}
%

\section{Preliminaries}

Choosing a new system of local coordinates leaving $x_0=\text{const.}$ to be invariant one can assume that
\[
p(x,\xi)=\xi_0^m+a_2(x,\xi')\xi_0^{m-2}+\cdots+a_m(x,\xi')
\]
and hence $\Sigma\subset \{\xi_0=0\}$. Thus near $\rho$ we may assume that $\Sigma$ is defined by $b_0(x,\xi)=\cdots =b_k(x,\xi)=0$ where $b_0=\xi_0$, $b_j=b_j(x,\xi')$, $1\leq j\leq k$ and $db_j$ are linearly independent at $\rho'$ where $\rho'$ stands for $({\bar x},{\bar \xi}')$ when $\rho=({\bar x},{\bar \xi})$. Recall that the localization $p_{\rho}(x,\xi)$ is a homogeneous hyperbolic polynomial of degree $m$ in $(x,\xi)$ in the direction  $(0,\theta)\in \R^{n+1}\times \R^{n+1}$. 
\begin{lemma}[{\cite[Lemma 1.1.3]{KN}}]
\label{lem:KNsp}
The next two conditions are equivalent.
\begin{description}
\item{\rm(i)} $C_{\rho}\cap T_{\rho}\Sigma=\{0\}$,
\item{\rm(ii)} $\Gamma_{\rho}\cap (T_{\rho}\Sigma)^{\sigma}\cap \lr{(0,\theta)}^{\sigma}\neq\emptyset$.
\end{description}
\end{lemma}
Assume $C_{\rho}\cap T_{\rho}\Sigma=\{0\}$ then thanks to Lemma \ref{lem:KNsp} there exists $0\neq X\in \Gamma_{\rho}\cap (T_{\rho}\Sigma)^{\sigma}\cap \lr{(0,\theta)}^{\sigma}$. Since $(T_{\rho}\Sigma)^{\sigma}$ is spanned by $H_{b_j}(\rho)$, $j=0,\ldots,k$ one can write
\begin{equation}
\label{eq:jikansu}
X=\sum_{j=0}^k\al_jH_{b_j}(\rho)
\end{equation}
where $\al_0=0$ because $X\in \lr{(0,\theta)}^{\sigma}$. This proves $\dif_{x_0}b_j(\rho')\neq 0$ with some $1\leq j\leq k$. Indeed if not we would have $X=(x,0,\xi')$ while denoting $p_{\rho}(x,\xi)=\prod_{j=1}^m(\xi_0-\Lambda_j(x,\xi'))$ we see $\Gamma_{\rho}=\{(x,\xi)\mid \xi_0>\max_j{\Lambda_j(x,\xi')}\}$ (for example \cite[Lemma 8.7.3]{Hobook:iti}) and we would have $\Lambda_j(x,\xi')<0$ which contradicts $\sum_{j=1}^m\Lambda_j(x,\xi')=0$. Renumbering, if necessary, one can assume $\dif_{x_0}b_1(\rho')\neq 0$ so that
\[
b_1(x,\xi')=(x_0-f_1(x',\xi'))e_1(x,\xi'),\quad e_1(x,\xi')\neq 0.
\]
Writing $b_j(x,\xi')=b_j(f_1(x',\xi'),x_1,\ldots,x_n,\xi')+c_j(x,\xi')b_1(x,\xi')$ we may assume $b_j(x,\xi')$, $2\leq j\leq k$ are independent of $x_0$. Since 
$p(x,\xi)$ vanishes on $\Sigma$ of order $m$ one can write with $b=(b_0,b_1,\ldots,b_k)=(b_0,b')$
\begin{equation}
\label{eq:ptob}
p(x,\xi)=b_0^m+\sum_{|\al|=m, \al_0\leq m-2}{\tilde a}_{\al}(x,\xi')b(x,\xi)^{\al}.
\end{equation}
Let ${\hat b}_j$ be  defined by $b_j(\rho+\mu X)=\mu {\hat b}_j(X)+O(\mu^2)$  and with ${\hat b}=(\xi_0,{\hat b}_1,\ldots,{\hat b}_k)$ we have $p_{\rho}(X)=q({\hat b}(X))$ where 
\[
q(\zeta)=\zeta_0^m+\sum_{|\al|=m, \al_0\leq m-2}{\tilde a}_{\al}(\rho')\zeta^{\al},\quad \zeta=(\zeta_0,\zeta_1,\ldots,\zeta_k)=(\zeta_0,\zeta')
\]
is a strictly hyperbolic polynomial in the direction $(1,0,\ldots,0)\in\R^{k+1}$ by \eqref{eq:assump}. Denote ${\tilde q}(\zeta; x,\xi')=q(\zeta)+\sum a_{\al}(x,\xi')\zeta^{\al}$
with $a_{\al}(x,\xi')={\tilde a}_{\al}(x,\xi')-{\tilde a}(\rho')$ 
and hence we have $p(x,\xi)={\tilde q}(b(x,\xi);x,\xi')$. 
\begin{lemma}
\label{lem:blow:a} There are $m$ real valued functions $\lambda_1(x,\xi')\leq \lambda_2(x,\xi')\leq \cdots\leq \lambda_m(x,\xi')$ defined in a conic neighborhood of $\rho'$ such that  
\begin{align*}
&p(x,\xi)=\prod_{j=1}^m\big(\xi_0-\lambda_j(x,\xi')\big),\quad|\lambda_j(x,\xi')|\leq C|b'(x,\xi')|,\\
&|\lambda_i(x,\xi')-\lambda_j(x,\xi')|\geq c\, |b'(x,\xi')|,\;\;(i\neq j)
\end{align*}
with some $c>0$, $C>0$.
\end{lemma}
\noindent
Proof: The first assertion is clear because $p(x,\xi)$ is a hyperbolic polynomial in the direction $\xi_0$. Note that ${\tilde q}(\zeta;\rho')=0$ has $m$ real distinct roots for  $\zeta'\neq 0$ then by Rouch\'e's theorem ${\tilde q}(\zeta_0,\zeta';x,\xi')=0$ has $m$ real distinct roots $\zeta_0=\lambda_j(\zeta';x,\xi')$ if $|\xi'-\rho'|$ is sufficiently small which are of homogeneous of degree $1$ in $\zeta'$ and $0$ in $\xi'$. It is easy to check that $|\lambda_j(\zeta';x,\xi')|\leq C|\zeta'|$ and $|\lambda_i(\zeta';x,\xi')-\lambda_j(\zeta';x,\xi')|\geq c\,|\zeta'|$ $(i\neq j)$ with some $c>0, C>0$. Since $\{|\zeta'|=1\}$ is compact we end the proof.
\qed
%

\section{Basic weights (energy estimates) }

We first introduce  symbol classes of pseudodifferential operators which will be used in this paper. Denote $\xim^2=\gamma^2+|\xi|^2$ where $\gamma\geq 1$ is a positive parameter.
\begin{definition}
\label{dfn:gmarus}\rm Let $W=W(x,\xi;\gamma)>0$ be a positive function. We define $S_{(s)}(W,g_0)$ to be the set of all $a(x,\xi;\gamma)\in C^{\infty}(\R^{n+1}\times \R^{n+1})$ such that
\begin{equation}
\label{eq:marus}
|\dif_x^{\be}\dif_{\xi}^{\al}a(x,\xi;\ga)|\leq CA^{|\al+\be|}|\al+\be|!^sW\xim^{-|\al|}
\end{equation}
and $S_{\lr1}(W,\bg)$ to be the set of all $a(x,\xi;\gamma)$ such that 
we have
\begin{equation}
\label{eq:bosai}
|\dif_x^{\be}\dif_{\xi}^{\al}a|\leq CA^{|\al+\be|}W(|\al+\be|+|\al+\be|^s\xim^{-\delta/2})^{|\al+\be|}\xim^{-\rho|\al|+\delta|\be|}
\end{equation}
for any $\al, \be\in\N^{n+1}$ with positive constants $C, A>0$ independent of $\ga\geq 1$. If $a(x,\xi;\ga)$ satisfies \eqref{eq:marus} {\rm (resp.\eqref{eq:bosai})} in a conic open set $U\subset \R^{n+1}\times (\R^{n+1}\setminus\{0\})$ we say $a(x,\xi;\ga)\in S_{(s)}(W,g_0)$ {\rm (resp. $S_{\lr1}(W,\bg)$)} in $U$. We often write $a(x,\xi)$ for $a(x,\xi;\ga)$  dropping $\ga$. 
\end{definition}
Note that $g_0$ and ${\bar g}$ is the metric defining the symbol class $S_{1,0}$ and $S_{\rho,\delta}$ respectively. It is clear that  one may replace 
$(|\al+\be|+|\al+\be|^s\xim^{-\delta/2})^{|\al+\be|}$ by $
|\al+\be|!(1+|\al+\be|^{s-1}\xim^{-\delta/2})^{|\al+\be|}$ in \eqref{eq:bosai}, still defining the same symbol class.

Since $p(x,\xi)$ is a polynomial in $\xi$ of degree $m$ it is clear that $p(x,\xi)\in S_{(s)}(\xim^m,g_0)$. Since $b_j(x,\xi')$ are defined only in a conic neighborhood of $\rho'=({\bar x}, {\bar\xi}')$ we extend such symbols to $\R^{n+1}\times\R^{n+1}$. Let 
$\chi(t)\in \gamma^{(s)}(\R)$ be $1$ for $|t|<c/2$ and $0$ for $|t|>c$ with small   $0<c<1/2$ and set
\[
\left\{\begin{array}{ll}
y(x)=\chi(|x-{\bar x}|)(x-{\bar x})+{\bar x},\\
\eta'(\xi)=\chi(|\xi'\xim^{-1}-{\bar\xi}'|)(\xi'-\xim{\bar\xi}')+\xim{\bar\xi}'.
\end{array}\right.
\]
Then it is easy to see $\eta'$, $b_j(y,\eta')\in S_{(s)}(\xim,g_0)$ and $\xim/C\leq |\eta'|\leq C\xim$ with some $C>0$.  In what follows we denote $b_j(y,\eta')$ by $b_j(x,\xi)$.

We now define $w(x,\xi)$, $\omega(x,\xi)$ by
\[
\left\{\begin{array}{ll}w(x,\xi)=\big(\sum_{j=1}^kb_j(x,\xi)^2\xim^{-2}
+\xim^{-2\delta}\big)^{1/2},\\[4pt]
\omega=\big(\phi^2+\xim^{-2\delta}\big)^{1/2},\;\;\phi=\sum_{j=1}^k\alpha_jb_j(x,\xi)\xim^{-1}.
\end{array}\right.
\]
Here we recall \eqref{eq:jikansu}, that is
\begin{equation}
\label{eq:kuri}
H_{\phi}(\rho)\in \Gamma_{\rho}.
\end{equation}
In what follows we assume that $0<\delta<\rho<1$ verifies $\rho+\delta=1$ and $1<s$ satisfies $0<s-1< (1-\rho)/2\rho$.
\begin{lemma}
\label{lem:www}There exist $C, A>0$ such that 
\[
|\dif_{x}^{\be}\dif_{\xi}^{\al}w|\leq CA^{|\al+\be|}(|\al+\be|+|\al+\be|^s\xim^{-\delta/2})^{|\al+\be|}w\xim^{-\rho|\al|+\delta|\be|}
\]
that is $w\in S_{\lr{1}}(w,\bg)$. We also have $\omega^{\pm 1}\in S_{\lr{1}}(\omega^{\pm 1},\bg)$.
\end{lemma}
We first remark an easy lemma.
\begin{lemma}
\label{lem:aa}Let $M>0$ be such that $2\big(1+4\sum_{j=0}^{\infty}(j+1)^{-2}\big)M\leq 1/2$ and $\Gamma_1(k)=Mk!/k^{3}$, $k\in {\mathbb N}$ where $\Gamma(0)=M$. Then we have
\[
\sum_{\al'+\al''=\al}\binom{\al}{\al'}\Gamma_1(|\al'|)\Gamma_1(|\al''|)\leq \Gamma_1(|\al|)/2.
\]
\end{lemma}
%
\noindent
Proof of Lemma \ref{lem:www}: It suffices to prove the assertion for $\epsilon w$ with small $\epsilon>0$ so that one can assume $|w|\leq 1$. Thus  with $w^2=F$  there is $A_1>0$ such that 
\[
|\dif_{x}^{\al}\dif_{\xi}^{\be}F|\leq A_1^{|\al+\be|}\Gamma_1(|\al+\be|)|\al+\be|^{(s-1)|\al+\be|}\xim^{-|\be|}
\]
holds for any $\al, \be$. Noting $|\dif_{x}^{\al}\dif_{\xi}^{\be}w|\leq C_{\al\be}w^{1-|\al+\be|}\xim^{-|\be|}$ for any $\al,\be$ we choose $A\geq 2A_1$ so that $
C_{\al\be}\leq A^{|\al+\be|}\Gamma_1(|\al+\be|)$ for $|\al+\be|\leq 4$ then we have
\begin{equation}
\label{eq:kinou}
\begin{split}
|\dif_{x}^{\al}\dif_{\xi}^{\be}w|
\leq A^{|\al+\be|}\Gamma_1(|\al+\be|)w\xim^{-\rho|\be|+\delta|\al|}\\
\times (w^{-1}\xim^{-\delta}+|\al+\be|^{s-1}\xim^{-\delta/2})^{|\al+\be|}.
\end{split}
\end{equation}
Suppose that \eqref{eq:kinou} holds for $|\al+\be|\leq k$, $4\leq k$ and let $|\al+\be|=k+1\geq 4$. Noting 
\[
2w\dif_{x}^{\al}\dif_{\xi}^{\be}w=-\sum_{1\leq |\al'+\be'|\leq k} \binom{\al}{\al'}\binom{\be}{\be'}\dif_{x}^{\al'}\dif_{\xi}^{\be'}w\dif_{x}^{\al-\al'}\dif_{\xi}^{\be-\be'}w+\dif_{x}^{\al}\dif_{\xi}^{\be}F
\]
and $w^{-1}\geq 1$, applying Lemma \ref{lem:aa} we see that $w|\dif_{x}^{\al}\dif_{\xi}^{\be}w|$ is bounded by
\begin{align*}
\frac{1}{2}A^{|\al+\be|}\Gamma_1(|\al+\be|)w^2\xim^{-\rho|\be|+\delta|\al|}(w^{-1}\xim^{-\delta}+|\al+\be|^{s-1}\xim^{-\delta/2})^{|\al+\be|}\\
+A_1^{|\al+\be|}
\Gamma_1(|\al+\be|)w^{2}|\al+\be|^{(s-1)|\al+\be|}(w^{-2}\xim^{-\delta|\al+\be|})\xim^{-\rho|\be|+\delta|\al|}.
\end{align*}
Since we have $w^{-2}\xim^{-\delta|\al+\be|}\leq \xim^{-\delta(|\al+\be|-2)}\leq \xim^{-\delta|\al+\be|/2}$ 
if $|\al+\be|\geq 4$ then taking 
$A^{|\al+\be|}/2+A_1^{|\al+\be|}\leq A^{|\al+\be|}$ into account we conclude that \eqref{eq:kinou} holds for $|\al+\be|=k+1$. Therefore noting $w^{-1}\xim^{-\delta}\leq 1$ we get
\[
|\dif_{x}^{\al}\dif_{\xi}^{\be}w|\leq A^{|\al+\be|}\Gamma_1(|\al+\be|)w\xim^{-\rho|\be|+\delta|\al|}
 (1+|\al+\be|^{s-1}\xim^{-\delta/2})^{|\al+\be|}.
\]
The assertion for $\omega$ is proved similarly. As for $\omega^{-1}$, using
\begin{align*}
\omega|\dif_{x}^{\al}\dif_{\xi}^{\be}\omega^{-1}|\leq \sum\binom{\al}{\al'}\binom{\be}{\be'}A^{|\al+\be|}\Gamma_1(|\al'+\be'|)\Gamma_1(|\al+\be|-|\al'+\be'|)\\
\times \xim^{-\rho|\be|+\delta|\al|}(w^{-1}\xim^{-\delta}+|\al+\be|^{s-1}\xim^{-\delta/2})^{|\al+\be|}
\end{align*}
the proof follows from induction on $|\al+\be|$.
\qed

\medskip
We now introduce a basic weight symbol which plays a key role in obtaining energy estimates:
\[
\psi=\xim^{\kappa}\log{(\phi+\omega)},\quad \kappa=\rho-\delta.
\]
\begin{lemma}
\label{lem:logphi}We have $(\phi+\omega)^{\pm 1}\in S_{\lr{1}}((\phi+\omega)^{\pm 1},\bg)$. We have also $\psi\in S_{\lr{1}}(\xim^{\kappa}\log{\xim},{\bar g})$. Moreover $\dif_x^{\be}\dif_{\xi}^{\al}
\psi\in S_{\lr{1 }}(\omega^{-1}\xim^{\kappa-|\al|},{\bar g})$ for $|\al+\be|=1$.
\end{lemma}
\noindent
Proof: With $W=\phi+\omega$ we put for $|\al+\be|=1$
\begin{equation}
\label{eq:kankei}
\dif_x^{\be}\dif_{\xi}^{\al}W=\frac{\dif_x^{\be}\dif_{\xi}^{\al}\phi}{\omega}W+\frac{\dif_x^{\be}\dif_{\xi}^{\al}\xim^{-2\delta}}{2\omega}=\Phi^{\al}_{\be}W+\Psi^{\al}_{\be}.
\end{equation}
We examine $\dif_x^{\be}\dif_{\xi}^{\al}\phi\in S_{\lr{1}}(\omega\xim^{-\rho|\al|+\delta|\be|},{\bar g})$ for $|\al+\be|=1$. Indeed noting $\omega^{-1}\xim^{-\delta}\leq 1$ we have
\begin{align*}
|\dif_x^{\nu+\be}\dif_{\xi}^{\mu+\al}\phi|
\leq CA^{|\mu+\nu|}\omega|\mu+\nu|^{s|\mu+\nu|}\xim^{-\delta|\mu+\nu|}\xim^{-\rho|\mu+\al|+\delta|\nu+\be|}\\
\leq CA^{|\mu+\nu|}\omega \xim^{-\rho|\al|+\delta|\be|}(|\mu+\nu|^s\xim^{-\delta/2})^{|\mu+\nu|}\xim^{-\rho|\mu|+\delta|\nu|}.
\end{align*}
Since $\omega^{-1}\in S_{\lr{1}}(\omega^{-1},\bg)$ one can find $A_1>0$ such that
\begin{align*}
|\dif_x^{\nu}\dif_{\xi}^{\mu}\Phi^{\al}_{\be}|\leq A_1^{|\mu+\nu|+1}\xim^{-\rho|\al+\mu|+\delta|\be+\nu|}|\mu+\nu|!\\
\times (1+|\mu+\nu|^{s-1}\xim^{-\delta/2})^{|\mu+\nu|},\;\;\forall \mu,\nu
\end{align*}
holds for $|\al+\be|=1$. Since $\xim^{-2\delta}\leq W$ similar arguments prove
\begin{equation}
\label{eq:Psino}
\begin{split}
|\dif_x^{\nu}\dif_{\xi}^{\mu}\Psi^{\al}_{\be}|
\leq A_1^{|\mu+\nu|+1}W\xim^{-\rho|\al+\mu|+\delta|\be+\nu|}|\mu+\nu|!\\
\times (1+|\mu+\nu|^{s-1}\xim^{-\delta/2})^{\mu+\nu|},\;\;\forall \mu,\nu.
\end{split}
\end{equation}
Now suppose 
\begin{equation}
\label{eq:Wno}
\begin{split}
|\dif_x^{\be}\dif_{\xi}^{\al}W|\leq CA_2^{|\al+\be|}W\xim^{-\rho|\al|+\delta|\be|}|\al+\be|!\\
\times (1+|\al+\be|^{s-1}\xim^{-\delta/2})^{|\al+\be|}
\end{split}
\end{equation}
holds for $|\al+\be|\leq \ell$ and letting $|\al+\be+e_1+e_2|=\ell+1$ we see
\begin{align*}
|\dif_{x}^{\be+e_2}\dif_{\xi}^{\al+e_1}W|
\leq C\sum\binom{\al}{\al'}\binom{\be}{\be'}A_1^{|\al-\al'+\be-\be'|+1}A_2^{|\al'+\be'|}\\
\times |\al-\al'+\be-\be'|!|\al'+\be'|!W\xim^{-\rho|\al+e_1|+\delta|\be+e_2|}\\
\times (1+|\al+\be|^{s-1}\xim^{-\delta/2})^{|\al+\be|}\\
+A_1^{|\al+\be|+1}|\al+\be|!W(1+|\al+\be|^{s-1}\xim^{-\delta/2})^{|\al+\be|}\xim^{-\rho|\al+e_1|+\delta|\be+e_2|}\\
\leq \big(CA_2^{|\al+\be|+1}A_1(A_2-A_1)^{-1}+A_1^{|\al+\be|+1}\big)\\
\times |\al+\be|!W(1+|\al+\be|^{s-1}\xim^{-\delta/2})^{|\al+\be|}\xim^{-\rho|\al+e_1|+\delta|\be+e_2|}.
\end{align*}
Thus it suffices to choose $A_2$ so that $A_1(A_2-A_1)^{-1}+C^{-1}(A_1A_2^{-1})\leq 1$ to conclude $\phi+\omega\in S_{\lr{1}}(\phi+\omega,\bg)$. As for $(\phi+\omega)^{-1}$ it suffices to repeat the  proof of Lemma \ref{lem:www}.

We turn to the next assertion. From $\xim^{-2\delta}/C\leq \phi+\omega\leq C$ it is clear $|\psi|\leq \xim^{\kappa}\log{\xim}$. Since $\dif_x^{\be}\dif_{\xi}^{\al}\log{(\phi+\omega)}=\dif_x^{\be}\dif_{\xi}^{\al}(\phi+\omega)/(\phi+\omega)$ for $|\al+\be|=1$ and $(\phi+\omega)^{\pm 1}\in S_{\lr{1}}((\phi+\omega)^{\pm 1},{\bar g})$ we see $\psi\in S_{\lr{1}}(\xim^{\kappa}\log{\xim}, {\bar g})$. Since $\dif_x^{\be}\dif_{\xi}^{\al}\phi\in S_{\lr{1}}(\xim^{-|\al|},\bg)$ for $|\al+\be|=1$ and $\omega^{-1}\in S_{\lr{1}}(\omega^{-1},\bg)$  it follows from \eqref{eq:kankei} and \eqref{eq:Wno} that
\begin{align*}
|\dif_x^{\nu}\dif_{\xi}^{\mu}(\dif_x^{\be}\dif_{\xi}^{\al}W)|\leq CA^{|\mu+\nu|}\omega^{-1}W\xim^{-|\al|}|\mu+\nu|!\\
\times (1+|\mu+\nu|^{s-1}\xim^{-\delta/2})^{\mu+\nu|}\xim^{-\rho|\mu|+\delta|\nu|}
\end{align*}
which proves the second assertion.
\qed
%

\section{Composition formula (energy estimates)}

In studying ${\rm Op}(e^{\psi})P{\rm Op}(e^{-\psi})={\rm Op}(e^{\psi}\#P\#e^{-\psi})$, if $\psi\in S^{\kappa'}_{\rho,\delta}$ with $\kappa'<\rho-\delta$ then one can apply the calculus obtained in \cite{NTa} to get an asymptotic formula of $e^{\psi}\#P\#e^{-\psi}$, where the proof is based on the almost analytic extension of symbols and the Stokes' formula using a space $\rho-\delta-\kappa'>0$. In the present  case $\psi\in S^{\kappa}_{\rho,\delta}$ there is no space  between $\kappa$ and $\rho-\delta$ and then, introducing a small parameter $\epsilon>0$, we carefully estimate $e^{\epsilon\psi}\#p\#e^{-\epsilon\psi}$ directly  to obtain the composition formula in Theorem \ref{thm:matome:6} below. 

We denote $a(x,\xi;\gamma,\epsilon)\in \epsilon^{\kappa_1}S_{\lr{1}}(W,\bg)$ if $\epsilon^{-\kappa_1}a\in S_{\lr{1}}(W,\bg)$ uniformly in $0<\epsilon\ll 1$. Our aim in this section is to give a sketch of the proof of
\begin{theorem}
\label{thm:matome:6}Let $p(x,\xi)\in S_{(s)}(\xim^m,g_0)$. Then there exists $\epsilon_0>0$ such that one can find $K=1+r$, $r\in \sqrt{\epsilon}\, S(1,{\bar g})$ and 
$\gamma_0(\epsilon)>0$ for $0<\epsilon\leq \epsilon_0$  so that  we have
\begin{align*}
e^{\epsilon\psi}\#p\#e^{-\epsilon\psi}\#K=\sum_{|\al+\be|\leq m}\frac{\epsilon^{|\al+\be|}}{\al!\be!}p^{(\be)}_{(\al)}(-i\nabla_{\xi}\psi)^{\al}(i\nabla_{x}\psi)^{\be}\\
+\sum_{1\leq |\al+\be|\leq m}\epsilon^{|\al+\be|+1/2}\, c^{\al}_{\be}p^{(\be)}_{(\al)}+R
\end{align*}
where $p^{(\be)}_{(\al)}=\dif_{\xi}^{\be}\dif_x^{\al}p$ and $c^{\al}_{\be}\in S(\xim^{\rho|\be|-\delta|\al|},{\bar g})$,  $R\in S(\xim^{m-\delta(m+1)},{\bar g})$. Moreover we have $c^{\al}_{\be}\in S(\omega^{-1}\xim^{\kappa-|\al|},\bg)$ for $|\al+\be|=1$. In particular $e^{\epsilon\psi}\#p\#e^{-\epsilon\psi}\in S(\xim^m,\bg)$.
\end{theorem}
%

\subsection{Estimates of symbol $e^{\epsilon\psi}$}

Let $H=H(x,\xi;\gamma)>0$ be a positive function. Assume that $f$ satisfies 
\begin{align*}
&\big|\dif_x^{\nu}\dif_{\xi}^{\mu}f\big|\leq  C_0A_0^{|\mu+\nu|}(|\mu+\nu|-1)!\\
&\times (1+(|\mu+\nu|-1)^{s-1}\xim^{-\delta/2})^{|\mu+\nu|-1}H\xim^{\delta|\nu|-\rho|\mu|}
\end{align*}
for $|\mu+\nu|\geq 1$. Set $\Omega^{\al}_{\be}=e^{-f}\dif_x^{\be}\dif_{\xi}^{\al}e^{f}$ then we have
\begin{lemma}
\label{lem:bb1}  Notations being as above. There exist $A_i, C>0$ such that the following estimate holds for $|\al+\be|\geq 1$;
\begin{eqnarray*}
&\big|\dif_x^{\nu}\dif_{\xi}^{\mu}\Omega^{\al}_{\be}\big|\leq CA_1^{|\nu+\mu|}A_2^{|\al+\be|}\xim^{\delta|\be+\nu|-\rho|\al+\mu|}\\
&\times \sum_{j=1}^{|\al+\be|}H^{|\al+\be|-j+1}(|\mu+\nu|+j)!(1+(|\mu+\nu|+j)^{s-1}\xim^{-\delta/2})^{|\mu+\nu|+j}.
\end{eqnarray*}
\end{lemma}
\newtheorem{cor}{Corollary}[section]
\begin{cor}
\label{cor:bb} We have with some $A, C>0$
\begin{align*}
|\dif_x^{\be}\dif_{\xi}^{\al}e^{f}|
\leq CA^{|\al+\be|}\xim^{\delta|\be|-\rho|\al|}
 \big(H+|\al+\be|+|\al+\be|^s\xim^{-\delta/2}\big)^{|\al+\be|}e^{f}.
\end{align*}
Moreover for $|\al+\be|\geq 1$, $|\dif_x^{\be}\dif_{\xi}^{\al}e^{f}|$ is bounded by
\begin{eqnarray*}
CHA^{|\al+\be|}\xim^{\delta|\be|-\rho|\al|}
 \big(H+|\al+\be|+|\al+\be|^s\xim^{-\delta/2}\big)^{|\al+\be|-1}e^{f}.
\end{eqnarray*}
\end{cor}
\begin{cor}
\label{cor:kaigo}
Notations being as above. We have for $|\al+\be|\geq 1$
\[
\Omega^{\al}_{\be}\in S_{\lr{1}}(H(H+|\al+\be|+|\al+\be|^s\xim^{-\delta/2})^{|\al+\be|-1}\xim^{-\rho|\al|+\delta|\be|},\bg).
\]
\end{cor}
\begin{cor}
\label{cor:rorei} Let $\omega^{\al}_{\be}=e^{-\epsilon\psi}\dif_x^{\be}\dif_{\xi}^{\al}e^{\epsilon\psi}$. Then there exists $\gamma_0(\epsilon)>0$ such that 
$\omega^{\al}_{\be}\in \epsilon^{|\al+\be|}S_{\lr{1}}(\xim^{\rho|\be|-\delta|\al|},\bg)$ for $\gamma\geq \gamma_0(\epsilon)$.
\end{cor}
%

\subsection{Estimates of $(b e^{\epsilon\psi})\# e^{-\epsilon\psi}$}

Let $\chi(r)\in \gamma^{(s)}(\R)$ be $1$ in $|r|\leq 1/4$ and $0$ outside $|r|\leq 1/2$.  
Let $b\in S_{\lr{1}}(\omega^t\xim^m,\bg)$ and consider  
\begin{align*}
(b e^{\epsilon\psi})\# e^{-\epsilon\psi}=\int e^{-2i(\eta z -y\zeta)}b(X+ Y)e^{\epsilon\psi(X+Y)-\epsilon\psi(X+Z)}dYdZ\\
=b(X)+\int e^{-2i(\eta z -y\zeta)}b(X+ Y)\big(e^{\epsilon\psi(X+Y)-\epsilon\psi(X+Z)}-1\big)dYdZ
\end{align*}
where $Y=(y,\eta)$, $Z=(z,\zeta)$. Denoting ${\hat\chi}=\chi(\lr{\eta}\xim^{-1})\chi(\lr{\zeta}\xim^{-1})$, ${\tilde \chi}=\chi(|y|/4)\chi(|z|/4)$ we write
\begin{align*}
\int e^{-2i(\eta z -y\zeta)}b(X+ Y)\big(e^{\epsilon\psi(X+Y)-\epsilon\psi(X+Z)}-1\big)\{{\tilde\chi}{\hat \chi}+(1-{\tilde\chi}){\hat \chi}\}dYdZ\\
+\int e^{-2i(\eta z -y\zeta)}b(X+ Y)\big(e^{\epsilon\psi(X+Y)-\epsilon\psi(X+Z)}-1\big)(1-{\hat\chi})dYdZ.
\end{align*}
After the change of variables $Z\to Z+Y$ the first integral turns to 
\begin{equation}
\label{eq:sitamati}
\int e^{-2i(\eta z -y\zeta)}b(X+ Y)(e^{\epsilon\psi(X+Y)-\epsilon\psi(X+Y+Z)}-1){\hat\chi}_0dYdZ
\end{equation}
where we have set ${\hat\chi}_0={\tilde\chi}(y,z)\chi(\lr{\eta}\xim^{-1})\chi(\lr{\eta+\zeta}\xim^{-1})$.
\begin{lemma}
\label{lem:sakura} Let $\Psi(X,Y,Z)=\psi(X+Y)-\psi(X+Y+Z)$ then  
on the support of ${\hat \chi}_0$ one has  
\[
|\Psi(X,Y,Z)|
\leq C\xim^{\kappa}{\bar g}_X^{1/2}(Z).
\]
 On the support of ${\hat\chi}_0$  we have for $|\al+\be|\geq 1$
\begin{align*}
&|\dif_{(x,y)}^{\be}\dif_{(\xi,\eta)}^{\al}e^{\epsilon\Psi}|\leq \epsilon\,CA^{|\al+\be|}
\xim^{-\rho|\al|+\delta|\be|}\xim^{\kappa}{\bar g}_X^{1/2}(Z)\\
&\times\big(\epsilon\xim^{\kappa}{\bar g}_X^{1/2}(Z)+|\al+\be|+|\al+\be|^s\xim^{-\delta/2}\big)^{|\al+\be|-1}e^{\epsilon\Psi}.
\end{align*}
\end{lemma}
\noindent
Proof: The assertions follow from Lemma \ref{lem:honeone} and Corollary \ref{cor:kaigo}.
\qed

\medskip
Introducing the following differential operators and symbols
\[
\left\{\begin{array}{ll}
L=1+4^{-1}\xim^{2\rho}|D_{\eta}|^2+4^{-1}\xim^{-2\delta}|D_y|^2,\\[3pt]
M=1+4^{-1}\xim^{2\delta}|D_{\zeta}|^2+4^{-1}\xim^{-2\rho}|D_z|^2,\\[3pt]
\Phi=1+\xim^{2\rho}|z|^2+\xim^{-2\delta}|\zeta|^2
=1+\xim^{2\kappa}{\bar g}_X(Z),\\[3pt]
\Theta=1+\xim^{2\delta}|y|^2+\xim^{-2\rho}|\eta|^2=1+{\bar g}_X(Y)\end{array}\right.
\]
so that $\Phi^{-N}L^Ne^{-2i(\eta z-y\zeta)}=e^{-2i(\eta z-y\zeta)}$, $\Theta^{-\ell}M^{\ell}e^{-2i(\eta z-y\zeta)}=e^{-2i(\eta z-y\zeta)}$ we make integration by parts in \eqref{eq:sitamati}. Let $
F=b(X+Y)(e^{\epsilon\Psi}-1)$, $\chi^*=\chi(\epsilon\, \Phi)$, $\chi_*=1-\chi^*$ and note $|(\xim^{-\rho}\dif_{z})^{\al}(\xim^{\delta}\dif_{\zeta})^{\be}\chi^*|\leq C_{\al\be}$ with $C_{\al\be}$ independent of $\epsilon$. Consider 
\begin{equation}
\label{eq:hiseki}
\begin{split}
&\int e^{-2i(\eta z-y\zeta)}\chi^*\dif_x^{\be}\dif_{\xi}^{\al}F{\hat \chi}_0dYdZ\\
&=\int e^{-2i(\eta z-y\zeta)}\Phi^{-N}L^N\Theta^{-\ell}M^{\ell}(\chi^*\dif_x^{\be}\dif_{\xi}^{\al}F{\hat \chi}_0)dYdZ.
\end{split}
\end{equation}
Applying Corollary \ref{cor:bb} we can estimate the integrand of the right-hand side of \eqref{eq:hiseki};
\begin{align*}
|\Phi^{-N}L^N
\Theta^{-\ell}M^{\ell}(\chi^*\dif_x^{\be}\dif_{\xi}^{\al}F{\hat\chi}_0)|
\leq C_{\ell}A^{2N+|\al+\be|+\ell}\Phi^{-N}\Theta^{-\ell}\big\{\epsilon\xim^{\kappa}g_X^{1/2}(Z)\\
\times \big(\epsilon\xim^{\kappa}g_X^{1/2}(Z)+2N+|\al+\be|+(2N+|\al+\be|)^s\xim^{-\delta/2})^{2N+|\al+\be|-1}e^{\epsilon\Psi}\\
+|e^{\epsilon\Psi}-1|\big(2N+|\al+\be|+(2N+|\al+\be|^s\xim^{-\delta/2}\big)^{2N+|\al+\be|}\big\}\\
\times \omega^t(X+Y)\xim^m.
\end{align*}
Here we remark the following easy lemma.
\begin{lemma}
\label{lem:kantan}
Let $A\geq 0$, $B\geq 0$. Then there exists $C>0$ independent of $n, m\in\N$, $A$, $B$ such that
\begin{align*}
(A+n+m+(n+m)^sB)^{n+m}
\leq C^{n+m}(A+n+n^sB)^n(A+m+m^sB)^m.
\end{align*}
\end{lemma}
Since $|e^{\epsilon \Psi}-1|\leq C|\epsilon\Psi|\leq C\epsilon \Phi^{1/2}\leq C\sqrt{\epsilon}$ 
on the support of $\chi^*$, the right-hand side  can be estimated by
\begin{align*}
C_{\ell}A_1^{2N+|\al+\be|}\Phi^{-N}\Theta^{-\ell}\big\{\epsilon\xim^{\kappa}\bg_X^{1/2}(Z)(\epsilon\xim^{\kappa}\bg_X^{1/2}(Z)+2N-1\\+(2N-1)^s\xim^{-\delta/2})^{2N-1}(\epsilon\xim^{\kappa}\bg_X^{1/2}(Z)+|\al+\be|\\
+|\al+\be|^s\xim^{-\delta/2})^{|\al+\be|-1}+\sqrt{\epsilon}(\epsilon\xim^{\kappa}\bg_X^{1/2}(Z)+2N
+(2N)^s\xim^{-\delta/2})^{2N}\\
\times (\epsilon\xim^{\kappa}\bg_X^{1/2}(Z)+|\al+\be|
+|\al+\be|^s\xim^{-\delta/2})^{|\al+\be|}\big\}\omega^t(X+Y)\xim^m
\end{align*}
where we remark $\omega^{\pm 1}(X+Y)\leq C\omega^{\pm 1}(X)(1+\bg_X^{1/2}(Y))$ on the support of ${\hat\chi}_0$ and hence we have
$\omega^t(X+Y)\leq C\omega^t(X)\Theta^{t'}$ with some $t'$. Noting
\begin{align*}
A_1^{2N}\Phi^{-N}\big(\epsilon\xim^{\kappa}\bg_X^{1/2}(Z)+2N+(2N)^s\xim^{-\delta/2}\big)^{2N}\\
=\Big(\epsilon\frac{A_1\xim^{\kappa}\bg_X^{1/2}(Z)}{\Phi^{1/2}}+\frac{2A_1N}{\Phi^{1/2}}+\frac{A_1(2N)^s\xim^{-\delta/2}}{\Phi^{1/2}}\Big)^{2N}
\end{align*}
and
\begin{align*}
A_1^{2N}\Phi^{-N}\epsilon
\xim^{\kappa}g_X^{1/2}(Z)\big(\epsilon\xim^{\kappa}\bg_X^{1/2}(Z)+2N-1+(2N-1)^s\xim^{-\delta/2}\big)^{2N-1}\\
\leq \epsilon A_1\Big(\epsilon\frac{A_1\xim^{\kappa}
\bg_X^{1/2}(Z)}{\Phi^{1/2}}+\frac{2A_1N}{\Phi^{1/2}}+\frac{A_1(2N)^s\xim^{-\delta/2}}{\Phi^{1/2}}\Big)^{2N-1}
\end{align*}
we choose $N=N(z,\zeta,\xi)$ so that $
2A_1N={\bar c}\,\Phi^{1/2}$ with a small ${\bar c}>0$. Then noting that $\Phi\leq 1+2\xim^{2\rho}|z|^2+2\xim^{-2\delta}|\zeta|^2\leq C\xim^{2\rho}$ on the support of ${\hat\chi}_0$ and therefore
$\Phi^{(s-1)/2}\xim^{-\delta/2}\leq C\xim^{-\epsilon_1}\leq C\ga^{-\epsilon_1} $ with some $\epsilon_1>0$ we have
\begin{align*}
\Big(\epsilon\frac{A_1
\xim^{\kappa}\bg_X^{1/2}(Z)}{\Phi^{1/2}}+\frac{2A_1N}{\Phi^{1/2}}+\frac{A_1(2N)^s\xim^{-\delta/2}}{\Phi^{1/2}}\Big)^{2N-1}\\
\leq \big(A_1\epsilon +{\bar c}+{\bar c}^s\Phi^{(s-1)/2}\xim^{-\delta/2}\big)^{2N-1}\leq Ce^{-c_1\Phi^{1/2}}
\end{align*}
choosing ${\bar c}$ small and $\gamma\geq \gamma_0(\epsilon)$ large.  
%
%
On the other hand since $\xim^{\kappa}\bg_X^{1/2}(Z)\leq \Phi^{1/2}$ it is clear
\begin{align*}
(\epsilon\xim^{\kappa}\bg_X^{1/2}(Z)+|\al+\be|
+|\al+\be|^s\xim^{-\delta/2})^{|\al+\be|}e^{-c\,\Phi^{1/2}}\\
\leq CA^{|\al+\be|}(|\al+\be|+|\al+\be|^s\xim^{-\delta/2})^{|\al+\be|}e^{-c'\Phi^{1/2}}.
\end{align*}
Set $\ell'=\ell-t'$. Then noting $e^{-c'\Phi^{1/2}}\leq C\Phi^{-\ell'}$ we have
\begin{equation}
\label{eq:itiban}
\begin{split}
|\Phi^{-N}L^N
\Theta^{-\ell}M^{\ell}(\chi^*\dif_x^{\be}\dif_{\xi}^{\al}F{\hat\chi}_0)|\\
\leq 
\sqrt{\epsilon}\, CA^{|\al+\be|}(|\al+\be|+|\al+\be|^s\xim^{-\delta/2})^{|\al+\be|}
 \\
 \times\omega^t(X)\xim^{m-\rho|\al|+\delta|\be|}\Theta^{-\ell'}\Phi^{-\ell'}.
 \end{split}
\end{equation}
Finally choosing $\ell>t'+(n+1)/2$ and recalling $\int \Theta^{-\ell'}\Phi^{-\ell'}dYdZ=C$  we conclude
\begin{align*}
&\Big|\int e^{-2i(\eta z-y\zeta)}\chi^*\dif_x^{\be}\dif_{\xi}^{\al}F{\hat \chi}_0dYdZ\Big|\\
&\leq \sqrt{\epsilon}\, CA^{|\al+\be|}(|\al+\be|+|\al+\be|^s\xim^{-\delta/2})^{|\al+\be|}\omega^t\xim^{m-\rho|\al|+\delta|\be|}.
\end{align*}
We next consider
\begin{align*}
\int e^{-2i(\eta z-y\zeta)}\chi_*\dif_x^{\be}\dif_{\xi}^{\al}F{\hat \chi}_0dYdZ
\\
=\int e^{-2i(\eta z-y\zeta)}L^N\Phi^{-N}M^{\ell}\Theta^{-\ell}\chi_*\dif_x^{\be}\dif_{\xi}^{\al}F{\hat \chi}_0dYdZ.
\end{align*}
Similar arguments obtaining \eqref{eq:itiban} show that
\begin{align*}
|\Phi^{-N}L^N
\Theta^{-\ell}M^{\ell}(\chi_*\dif_x^{\be}\dif_{\xi}^{\al}F{\hat\chi}_0)|
\leq  CA^{|\al+\be|}(|\al+\be|+|\al+\be|^s\xim^{-\delta/2})^{|\al+\be|}\\
\times \omega^t(X)\xim^{m-\rho|\al|+\delta|\be|}\Theta^{-\ell'}\Phi^{-\ell'}e^{-c\,\Phi^{1/2}}.
\end{align*}
Since $\Phi^{1/2}\geq \epsilon^{-1/2}$ on the support of $\chi_*$ we see $e^{-c\,\Phi^{1/2}}\leq e^{-c\,\epsilon^{-1/2}}\leq C\sqrt{\epsilon}$ and this proves
\begin{align*}
&\Big|\dif_x^{\be}\dif_{\xi}^{\al}\int e^{-2i(\eta z-y\zeta)}F{\hat \chi}_0dYdZ\Big|\\
&\leq \sqrt{\epsilon}\, CA^{|\al+\be|}(|\al+\be|+|\al+\be|^s\xim^{-\delta/2})^{|\al+\be|}\omega^t\xim^{m-\rho|\al|+\delta|\be|}.
\end{align*}
We then consider 
\[
\int e^{-2i(\eta z-y\zeta)}(|y|^2+|z|^2)^{-N}(|D_{\zeta}|^2+|D_{\eta}|^2)^N F{\hat \chi}_1dYdZ
\]
 where $F=b(X+Y)(e^{\epsilon\psi(X+Y)-\epsilon\psi(X+Z)}-1)$ and ${\hat\chi}_1=(1-{\tilde\chi}){\hat\chi}$. Let $\kappa<\kappa_1<\rho$ then since
$|\psi(X+Y)|+|\psi(X+Z)|$ is bounded by $C\xim^{\kappa_1}$  and $\lr{\xi+\eta}_{\ga}\approx \xim$, $\lr{\xi+\zeta}_{\ga}\approx \xim$  on the support of ${\hat\chi}_1$ thanks to Corollary \ref{cor:bb} it is not difficult to show
\begin{align*}
\big|(|D_{\zeta}|^2+|D_{\eta}|^2)^N\dif_x^{\be}\dif_{\xi}^{\al}F{\hat \chi}_1\big|
\leq CA^{2N+|\al+\be|}\omega^t(X+Y)\xim^{m-\rho|\al|+\delta|\be|}\\
\times (\xim^{\kappa_1}+|\al+\be|
+|\al+\be|^s\xim^{-\delta/2})^{|\al+\be|}\\
\times \xim^{-2\rho N} (\xim^{\kappa_1}+2N+(2N)^s\xim^{-\delta/2})^{2N}e^{c\xim^{\kappa_1}}.
\end{align*}
Choose $N=c_1 \xim^{\rho}$ with small $c_1>0$ so that
\begin{align*}
A^{2N}\xim^{-2\rho N}\big( \xim^{\kappa_1}+2N+(2N)^s\xim^{-\delta/2}\big)^{2N}
\end{align*}
is bounded by $Ce^{-c\xim^{\rho}}$ and $\xim^{\kappa_1|\al+\be|}e^{-c\xim^{\rho}}$ is bounded by $CA^{|\al+\be|} e^{-c'\xim^{\rho}}$. Then noting  $\omega^t(X+Y)\leq C\omega^t(X)\xim^{t'}$ and $e^{-c'\xim^{\rho}}\leq \sqrt{\epsilon}\,C\xim^{-2(n+1)-t'}$ for $\ga\geq \ga_0(\epsilon)$ and that $\xim^{-2(n+1)}\int (|y|^2+|z|^2)^{-N}{\hat\chi}_1dYdZ\leq C$ 
we conclude
\begin{lemma}
\label{lem:syubu} Let ${\hat\chi}=\chi(\lr{\eta}\xim^{-1})\chi(\lr{\zeta}\xim^{-1})$. Then we have for $\ga\geq \ga_0(\epsilon)$
\begin{align*}
\Big|\dif_x^{\be}\dif_{\xi}^{\al}\int e^{-2i(z\eta-y\zeta)}b(X+Y)(e^{\epsilon\psi(X+Y)-\epsilon\psi(X+Z)}-1){\hat\chi}dYdZ\Big|\\
\leq \sqrt{\epsilon}\, CA^{|\al+\be|}(|\al+\be|+|\al+\be|^s\xim^{-\delta/2})^{|\al+\be|}\omega^t(X)\xim^{m-\rho|\al|+\delta|\be|}.
\end{align*}
\end{lemma}
Let us write 
\begin{align*}
1-{\hat\chi}=(1-\chi(\lr{\eta}\xim^{-1}))(1-\chi(\lr{\zeta}\xim^{-1}))
+(1-\chi(\lr{\eta}\xim^{-1}))\chi(\lr{\zeta}\xim^{-1})\\
+(1-\chi(\lr{\zeta}\xim^{-1}))\chi(\lr{\eta}\xim^{-1})
={\hat\chi}_2+{\hat\chi}_3+{\hat\chi}_4.
\end{align*}
Denoting $F=b(X+Y)(e^{\epsilon\psi(X+Y)-\epsilon\psi(X+Z)}-1)$ again we consider  
\begin{align*}
\int e^{-2i(\eta z-y\zeta)}\lr{\eta}^{-2N_2}\lr{\zeta}^{-2N_1}\lr{D_z}^{2N_2}\lr{D_y}^{2N_1}\\
\times \lr{y}^{-2\ell}\lr{z}^{-2\ell}\lr{D_{\zeta}}^{2\ell}\lr{D_{\eta}}^{2\ell}\dif_x^{\be}\dif_{\xi}^{\al}F{\hat\chi}_2\chi^*dYdZ
\end{align*}
where $\chi^*$ is either $\chi(\lr{\zeta}\lr{\eta}^{-1}/4)$ or $1-\chi(\lr{\zeta}\lr{\eta}^{-1}/4)$. If $\chi^*=\chi(\lr{\zeta}\lr{\eta}^{-1}/4)$ we choose $N_1=\ell$, $N_2=N$ and noting $\omega^t(X+Y)\leq C\lr{\eta}^{t'}$ with some $t'\geq 0$ and $|\psi(X+Y)|+|\psi(X+Z)|\leq C\lr{\eta}^{\kappa_1}$ with $\kappa<\kappa_1<\rho$ on the support of ${\hat \chi}_2\chi^*$ it is not difficult to see that 
\[
|\lr{\eta}^{-2N}\lr{\zeta}^{-2\ell}\lr{D_z}^{2N}\lr{D_y}^{2\ell}\lr{y}^{-2\ell}\lr{z}^{-2\ell}\lr{D_{\zeta}}^{2\ell}\lr{D_{\eta}}^{2\ell}\dif_x^{\be}\dif_{\xi}^{\al}F{\hat\chi}_2\chi^*|
\]
is bounded by
\begin{equation}
\label{eq:copi}
\begin{split}
C_{\ell}A^{2N+|\al+\be|}\lr{\eta}^{-2N}\lr{\zeta}^{-2\ell}\lr{y}^{-2\ell}\lr{z}^{-2\ell}\lr{\eta}^{m+t'+2\delta\ell}\lr{\eta}^{6\ell\rho}\\
\times (\lr{\eta}^{\kappa_1}+2N\lr{\eta}^{\delta}+N^s\lr{\eta}^{\delta/2})^{2N}\\
\times (\lr{\eta}^{\kappa_1}+\lr{\eta}^{\delta}|\al+\be|+\lr{\eta}^{\delta/2}|\al+\be|^s)^{|\al+\be|}
e^{C\lr{\eta}^{\kappa_1}}.
\end{split}
\end{equation}
Here writing 
\begin{align*}
A^{2N}\lr{\eta}^{-2N}\big( \lr{\eta}^{\kappa_1}+2N\lr{\eta}^{\delta}+N^s\lr{\eta}^{\delta/2}\big)^{2N}\\
=\Big(\frac{A \lr{\eta}^{\kappa_1}}{\lr{\eta}}+\frac{2AN}{\lr{\eta}^{\rho}}+\frac{AN^s\lr{\eta}^{\delta/2}}{\lr{\eta}}\Big)^{2N}
\end{align*}
we take $2N=c_1 \lr{\eta}^{\rho}$ with small $c_1>0$ so that the right-hand side is bounded by $Ce^{-c\lr{\eta}^{\rho}}$. Noting $\lr{\eta}^{\delta |\al+\be|}e^{-c\lr{\eta}^{\rho}}\leq CA_1^{|\al+\be|}|\al+\be|^{\delta |\al+\be|/\rho}e^{-c_1\lr{\eta}^{\rho}}$ and $\lr{\eta}^{\kappa_1 |\al+\be|}e^{-c\lr{\eta}^{\rho}}\leq CA_1^{|\al+\be|}|\al+\be|^{|\al+\be|}e^{-c_1\lr{\eta}^{\rho}}$ one sees that \eqref{eq:copi} is bounded by
\begin{align*}
C_{\ell}A_1^{|\al+\be|}\lr{\zeta}^{-2\ell}\lr{y}^{-2\ell}\lr{z}^{-2\ell}(|\al+\be|^{1+\delta/\rho}+|\al+\be|^{s+\delta/2\rho})^{|\al+\be|}e^{-c_1\lr{\eta}^{\rho}}.
\end{align*}
Similarly if $\chi^*=1-\chi(\lr{\zeta}\lr{\eta}^{-1}/4)$ choosing $N_1=N$, $N_2=\ell$ it is proved that \eqref{eq:copi} is estimated by 
\begin{align*}
C_{\ell}A_1^{|\al+\be|}\lr{\eta}^{-2\ell}\lr{y}^{-2\ell}\lr{z}^{-2\ell}(|\al+\be|^{1+\delta/\rho}+|\al+\be|^{s+\delta/2\rho})^{|\al+\be|}e^{-c_1\lr{\zeta}^{\rho}}.
\end{align*}
Thus taking $1+\delta/\rho=1/\rho$ and $s+\delta/2\rho\leq 1/\rho$ into account and recalling that $\xim\leq \lr{\eta}$, $\xim\leq  \lr{\zeta}$ on the support of ${\hat\chi}_2$ we get
\begin{lemma}
\label{lem:syu:b}We have
\begin{align*}
\Big|\dif_x^{\be}\dif_{\xi}^{\al}\int e^{-2i(z\eta-y\zeta)}b(X+Y)(e^{\epsilon\psi(X+Y)-\epsilon\psi(X+Z)}-1){\hat\chi}_2dYdZ\Big|\\
\leq  CA^{|\al+\be|}|\al+\be|^{|\al+\be|/\rho}e^{-c_1\xim^{\rho}}.
\end{align*}
\end{lemma}
Repeating similar arguments we can prove 
\begin{align*}
\Big|\dif_x^{\be}\dif_{\xi}^{\al}\int e^{-2i(\eta z -y\zeta)}b(X+Y)(e^{\epsilon\psi(X+Y)-\epsilon\psi(X+Z)}-1){\hat \chi}_idYdZ\Big|\\
\leq CA^{|\al+\be|}|\al+\be|^{|\al+\be|/\rho}e^{-c_1\xim^{\rho}}
\end{align*}
for $i=3,4$. We summarize what we have proved in
\begin{pro}
\label{pro:syu:AA} Let $b\in S_{\lr{1}}(\omega^t\xim^m,\bg)$ then we have
\[
(b e^{\epsilon\psi})\#e^{-\epsilon\psi}=b+\omega^t\,{\hat b}+R
\]
where ${\hat b}\in \sqrt{\epsilon}\, S_{\lr{1}}(\xim^m,\bg)$and $R\in S_{(1/\rho)}(e^{-c_1\xim^{\rho}},|dx|^2+|d\xi|^2)$, that is 
\[
|\dif_x^{\be}\dif_{\xi}^{\al}R|\leq CA^{|\al+\be|}|\al+\be|^{|\al+\be|/\rho}e^{-c_1\xim^{\rho}}.
\]
\end{pro}
%
\subsection{Proo of Theorem \ref{thm:matome:6}}

We start with the next lemma which is proved repeating similar arguments in the preceding subsection. 
\begin{lemma}
\label{lem:apsiomega} Let $a\in S_{(s)}(\xim^d,g_0)$ and $b\in S_{\lr{1}}(\xim^h,\bg)$. Then we have
\[
(be^{\epsilon\psi})\#a=\sum_{|\al+\be|<N}\frac{(-1)^{|\be|}}{(2i)^{|\al+\be|}\al!\be!}a^{(\be)}_{(\al)}(be^{\epsilon\psi})^{(\al)}_{(\be)}+b_Ne^{\epsilon\psi}+R
\]
where $b_N\in S_{\lr{1}}(\xim^{d+h-\delta N},{\bar g})$, $R\in S_{(1/\rho)}(e^{-c\xim^{\rho}}, |dx|^2+|d\xi|^2)$. For $a\#(be^{\epsilon\psi})$ similar assertion holds, where $(-1)^{|\be|}$ is replaced by $(-1)^{|\al|}$.
\end{lemma}
We can also prove
\begin{lemma}
\label{lem:kore} Let $R\in S_{(1/\rho)}(e^{-c\xim^{\rho}},|dx|^2+|d\xi|^2)$. Then we have
\begin{align*}
R\#e^{\pm\epsilon\psi},\;\;e^{\pm\epsilon\psi}\#R\in S_{(1/\rho)}(e^{-c'\xim^{\rho}}, |dx|^2+|d\xi|^2).
\end{align*}
\end{lemma}
\begin{cor}
\label{cor:sono} Let $R\in S_{(1/\rho)}(e^{-c\xim^{\rho}}, |dx|^2+|d\xi|^2)$. Then for any $t\in\R$ we have
\[
e^{\pm\epsilon\psi}\#R\#e^{\mp\epsilon\psi}\in S(\xim^t,g_0).
\]
\end{cor}
\begin{lemma}
\label{lem:matome:3} Let $p\in S_{(s)}(\xim^{d},g_0)$. Then one can write
\begin{align*}
(e^{\epsilon\psi})\# p=\sum_{|\al+\be|<N}\frac{(-1)^{|\be|}}{i^{|\al+\be|}\al!\be!}p^{(\be)}_{(\al)}\#(\omega^{\al}_{\be}e^{\epsilon\psi})+r_Ne^{\epsilon\psi}+R
\end{align*}
where $r_N\in S_{\lr{1}}(\xim^{d-\delta N},{\bar g})$, $R\in S_{(1/\rho)}(e^{-c\xim^{\rho}},|dx|^2+|d\xi|^2)$ and $\omega^{\al}_{\be}=e^{-\epsilon\psi}\dif_x^{\be}\dif_{\xi}^{\al}e^{\epsilon\psi}$.
\end{lemma}
\noindent
Proof: We first examine
\begin{equation}
\label{eq:masa}
p^{(\be)}_{(\al)}\omega
^{\al}_{\be}e^{\epsilon\psi}-\sum_{|\gamma+\delta|<N}\frac{(-1)^{|\gamma|}}{(2i)^{|\gamma+\delta|}\gamma!\delta!}p^{(\be+\gamma)}_{(\al+\delta)}\#(\omega^{\al+\delta}_{\be+\gamma}e^{\epsilon\psi})=r_{N,|\al+\be|}e^{\epsilon\psi}+R
\end{equation}
with $r_{N,|\al+\be|}\in S_{\lr{1}}(\xim^{d-\delta N},{\bar g})$. Indeed since $\dif_x^{\nu}\dif_{\xi}^{\mu}(\omega^{\al}_{\be}e^{\epsilon\psi})=\omega^{\al+\mu}_{\be+\nu}e^{\epsilon\psi}$ thanks to Lemma \ref{lem:apsiomega} one can write
\begin{align*}
&\sum_{|\gamma+\delta|<N}\frac{(-1)^{|\gamma|}}{(2i)^{|\gamma+\delta|}\gamma!\delta!}p^{(\be+\gamma)}_{(\al+\delta)}\#(\omega^{\al+\delta}_{\be+\gamma}e^{\epsilon\psi})
\\
&=\sum_{|\gamma'+\delta'|<2N}\frac{(-1)^{\gamma'|}}{(2i)^{|\gamma'+\delta'|}\gamma'!\delta'!}\big(\sum\binom{\gamma'}{\mu}\binom{\delta'}{\nu}(-1)^{|\mu+\nu|}\big)p^{(\be+\gamma')}_{(\al+\delta')}\omega^{\al+\delta'}_{\be+\gamma'}e^{\epsilon\psi}\\
&+r_{N,|\al+\be|}e^{\epsilon\psi}+R
\end{align*}
where $\sum\binom{\gamma'}{\mu}\binom{\delta'}{\nu}(-1)^{|\mu+\nu|}=0$ if $|\gamma'+\delta'|>0$ so that the right-hand side is 
\[
p^{(\be)}_{(\al)}\omega^{\al}_{\be}e^{\epsilon\psi}+r_{N,|\al+\be|}e^{\epsilon\psi}+R,\;\;r_{N,|\al+\be|}\in S_{\lr{1}}(\xim^{d-\delta N},{\bar g})
\]
which proves \eqref{eq:masa}. Now insert the expression of $p^{(\be)}_{(\al)}\omega^{\al}_{\be}e^{\epsilon\psi}$ in \eqref{eq:masa} into
\[
(e^{\epsilon\psi})\# p=\sum_{|\al+\be|<N}\frac{(-1)^{|\be|}}{(2i)^{|\al+\be|}\al!\be!}p^{(\be)}_{(\al)}\omega^{\al}_{\be}e^{\epsilon\psi}+r_Ne^{\epsilon\psi}+R
\]
which follows from Lemma \ref{lem:apsiomega} to get
\begin{align*}
\sum_{|\al'+\be'|<2N}\frac{(-1)^{|\be'|}}{(2i)^{|\al'+\be'|}\al'!\be'!}\big(\sum\binom{\al'}{\delta}\binom{\be'}{\gamma}\big)p^{(\be')}_{(\al')}\#(\omega^{\al'}_{\be'}e^{\epsilon\psi})+{\tilde r}_Ne^{\epsilon\psi}+R
\end{align*}
where ${\tilde r}_N\in  S_{\lr{1}}(\xim^{d-\delta N},\bg)$.  Here we note $\sum\binom{\al'}{\delta}\binom{\be'}{\gamma}=2^{|\al'+\be'|}$. It is clear that   $p^{(\be')}_{(\al')}\#(\omega^{\al'}_{\be'}e^{\epsilon\psi})=r'e^{\epsilon\psi}+R$ with $r'\in S_{\lr{1}}(\xim^{d-\delta N},\bg)$ for $|\al'+\be'|\geq N$ and hence we get the assertion.
\qed

\medskip
\noindent
Proof of Theorem \ref{thm:matome:6}: From Lemma \ref{lem:matome:3} we see
\begin{align*}
(e^{\epsilon\psi})\# p\#e^{-\epsilon\psi}=\sum_{|\al+\be|\leq m}\frac{(-1)^{|\be|}}{i^{|\al+\be|}\al!\be!}p^{(\be)}_{(\al)}\#((\omega^{\al}_{\be}e^{\epsilon\psi})\#e^{-\epsilon\psi})\\
+(r_me^{\epsilon\psi}+R)\#e^{-\epsilon \psi}
\end{align*}
where $(r_me^{\epsilon\psi}+R)\#e^{-\epsilon \psi}\in S(\xim^{m-\delta(m+1)},{\bar g})$ which follows from Propositions \ref{pro:syu:AA} and Lemma \ref{lem:kore}. Therefore Propositions \ref{pro:syu:AA} together with Corollary \ref{cor:rorei} gives 
\[
(e^{\epsilon\psi})\# p\#e^{-\epsilon\psi}=\sum_{|\al+\be|\leq m}\frac{(-1)^{|\be|}}{i^{|\al+\be|}\al!\be!}p^{(\be)}_{(\al)}\#(\omega^{\al}_{\be}+{\bar\omega}^{\al}_{\be})\\
+R
\]
where ${\bar \omega}^{\al}_{\be}\in \epsilon^{|\al+\be|+1/2} S_{\lr{1}}(\xim^{\rho|\be|-\delta|\al|},\bg)$, $R\in S(\xim^{m-\delta (m+1)},{\bar g})$ and 
\[
{\bar\omega}^{\al}_{\be}\in \epsilon^{3/2}S_{\lr{1}}(\omega^{-1}\xim^{\kappa-|\al|},{\bar g}),\quad |\al+\be|=1.
\]
 Since $e^{\epsilon\psi}\#e^{-\epsilon\psi}=1-r_1$, $r_1\in \sqrt{\epsilon}\, S(1,{\bar g})$ by Proposition \ref{pro:syu:AA} there exists $K=1+r$, $r\in\sqrt{\epsilon}\,S(1,{\bar g})$ such that $e^{\epsilon\psi}\#e^{-\epsilon\psi}\#K=1$ if $0<\epsilon\leq \epsilon_0$ is small (\cite[Theorem 3.2]{Bea} and  \cite[Theorem 2.6.27]{Ler} for example). Thus we have
\[
(e^{\epsilon\psi})\# p\#e^{-\epsilon\psi}\# K=p+\sum_{1\leq |\al+\be|\leq m}\frac{(-1)^{|\be|}}{i^{|\al+\be|}\al!\be!}p^{(\be)}_{(\al)}\#(\omega^{\al}_{\be}+{\bar\omega}^{\al}_{\be})\#K
+R.
\]
On the other hand it is clear $
(\omega^{\al}_{\be}+{\bar\omega}^{\al}_{\be})\#(1+r)=\omega^{\al}_{\be}+{\tilde \omega}^{\al}_{\be}$ with ${\tilde\omega}^{\al}_{\be}\in \epsilon^{|\al+\be|+1/2} S(\xim^{\rho|\be|-\delta|\al|},{\bar g})$. 
Since $\omega^{\al}_{\be}\in \epsilon S(\omega^{-1}\xim^{\kappa-|\al|},\bg)$ for $|\al+\be|=1$ it is also clear that ${\tilde\omega}^{\al}_{\be}\in \epsilon^{3/2} S(\omega^{-1}\xim^{\kappa-|\al|},{\bar g})$ for $|\al+\be|=1$. Note 
\begin{align*}
p^{(\be)}_{(\al)}\#(\omega^{\al}_{\be}+{\bar\omega}^{\al}_{\be})-\sum_{|\mu+\nu|\leq m-|\al+\be|}\frac{(-1)^{|\nu|}}{i^{|\mu+\nu|}\mu!\nu!}p^{(\be+\mu)}_{(\al+\nu)}(\omega^{\al}_{\be}+{\bar\omega}^{\al}_{\be})^{(\nu)}_{(\mu)}\\
\in S(\xim^{m-\delta|\al+\be|-\rho(m+1-|\al+\be|)},{\bar g})\subset S(\xim^{m-\delta (m+1)},\bg)
\end{align*}
and $(\omega^{\al}_{\be}+{\bar\omega}^{\al}_{\be})^{(\nu)}_{(\mu)}\in \gamma^{-\kappa|\mu+\nu|}S_{\lr{1}}(\xim^{\rho|\al+\mu|-\delta|\be+\nu|},\bg)$ which is contained in $ \epsilon^{|\al+\be|+|\mu+\nu|+1/2}S_{\lr{1}}(\xim^{\rho|\al+\mu|-\delta|\be+\nu|},\bg)$ if $|\mu+\nu|\geq 1$, $\gamma\geq \gamma_0(\epsilon)$  so that 
\[
(e^{\epsilon\psi})\# p\#e^{-\epsilon\psi}\# K=p+\sum_{1\leq |\al+\be|\leq m}\frac{(-1)^{|\be|}}{i^{|\al+\be|}\al!\be!}p^{(\be)}_{(\al)}(\omega^{\al}_{\be}+{\hat\omega}^{\al}_{\be})
+R
\]
where ${\hat\omega}^{\al}_{\be}\in \epsilon^{|\al+\be|+1/2} S(\xim^{\rho|\be|-\delta|\al|},{\bar g})$ and ${\hat\omega}^{\al}_{\be}\in \epsilon^{3/2} S(\omega^{-1}\xim^{\kappa-|\al|},{\bar g})$ for $|\al+\be|=1$. Now check $\omega^{\al}_{\be}$.
For $|\al+\be|=1$ we have $\omega^{\al}_{\be}=\epsilon(-i\nabla_{\xi}\psi)^{\al}(i\nabla_x\psi)^{\be}$. Let $|\al+\be|\geq 2$ then $\omega^{\al}_{\be}$ is a linear combination of terms $
(\epsilon\psi)^{(\al_1)}_{(\be_1)}\cdots (\epsilon \psi)^{(\al_s)}_{(\be_s)}$ 
with $\al_1+\cdots+\al_s=\al$, $\be_1+\cdots+\be_s=\be$, $|\al_i+\be_i|\geq 1$. If $|\al_i+\be_i|=1$ for all $i$ it is clear $\omega^{\al}_{\be}=\epsilon^{|\al+\be|}(-i\nabla_{\xi}\psi)^{\al}(i\nabla_x\psi)^{\be}$. If $|\al_j+\be_j|\geq 2$ for some $j$ so that $s\leq |\al+\be|-2$ then one has 
\[
(\epsilon\psi)^{(\al_1)}_{(\be_1)}\cdots (\epsilon\psi)^{(\al_s)}_{(\be_s)}\in S(\xim^{-\kappa+\rho|\be|-\delta|\al|}, {\bar g})\subset \gamma^{-\kappa} S(\xim^{\rho|\al|-\delta|\be|},\bg).
\]
Since we can assume  $\gamma^{-\kappa}\leq \epsilon^{|\al+\be|+1/2}$  for $\gamma\geq \gamma_0(\epsilon)$ we get the assertion.
\qed
%

\section{Energy estimates}

To obtain energy estimates we follow  \cite{KNW} where the main point is to derive microlocal energy estimates. We sketch how to get microlocal energy estimates. Let us denote 
\[
P_{\psi}={\rm Op}(e^{\epsilon\psi})P{\rm Op}(e^{-\epsilon\psi}){\rm Op}(K)
\]
of which principal symbol is given by $p_{\psi}=e^{\epsilon\psi}\#p\#e^{-\epsilon\psi}\#K$. In this section we say $a(x,\xi;\gamma,\epsilon)\in {\tilde S}(W,\bg)$ if $a\in S(W,\bg)$ for each fixed  $0<\epsilon\ll 1$. Let $a\in {\tilde S}(W,\bg)$ and let $N\in\N$ be given. Then with a fixed small $0<2\tau<\rho-\delta$ we have
\[
|\dif_x^{\be}\dif_{\xi}^{\al}a|\leq C_{\al\be}(\epsilon)W\xim^{-\rho|\al|+\delta|\be|}\leq C_{\al\be}\gamma^{-2\tau|\al+\be|}W\xim^{-(\rho-\tau)|\al|+(\delta+\tau)|\be|}
\]
where one can assume that $C_{\al\be}(\epsilon)\gamma^{-2\tau|\al+\be|}$ are arbitrarily small for $1\leq |\al+\be|\leq N$ taking $\ga$ large. 
\subsection{Symbol of $P_{\psi}$}

Define $h_j(x,\xi)$ by
\[
h_j(x,\xi)=\sum_{1\leq \ell_1<\ell_2<\cdots<\ell_j\leq m}|q_{\ell_1}|^2\cdots|q_{\ell_j}|^2,\quad q_j=\xi_0-\lambda_j(x,\xi').
\]
\begin{lemma}
\label{lem:preaa}
There exists $c>0$ such that
\[
h_{m-k}(x,\xi-i\epsilon\omega^{-1}\xim^{\kappa}\theta)\geq c(\epsilon\omega)^{2(j-k)}\xim^{2(j-k)}h_{m-j}(x,\xi-i\epsilon\omega^{-1}\xim^{\kappa}\theta)
\]
for $j=k,\ldots,m$ where $h_0=1$ and $1\leq k\leq m$.
\end{lemma}
\noindent
Proof: We show the case $k=1$. By definition $h_{m-1}(x,\xi-i\epsilon\omega^{-1}\xim^{\kappa}\theta)$ is bounded from below by 
\begin{align*}
2^{-1}(|q_i(x,\xi)|^2+|q_j(x,\xi)|^2+\epsilon^2\omega^{-2}\xim^{2\kappa})
\prod_{k\neq i,j}|q_{k}(x,\xi-i\epsilon\omega^{-1}\xim^{\kappa}\theta)|^2.
\end{align*}
From Lemma \ref{lem:blow:a} we have $|q_i(x,\xi)|^2+|q_j(x,\xi)|^2\geq c\, |b'(x,\xi)|^2$. Since
\begin{align*}
&|b'(x,\xi)|^2+\epsilon^2\omega^{-2}\xim^{2\kappa}
\\
&=\epsilon^2\omega^{-2}\xim^{2}(\epsilon^{-2}|b'(x,\xi)|^2\omega^2\xim^{-2}+\xim^{-4\delta})
\geq  c\,\epsilon^2\omega^2\xim^2
\end{align*}
with some $c>0$ because $C|b'(x,\xi)|^2\xim^{-2}\geq \phi^2$ and $\omega^2\geq \phi^2$ then it is clear that $h_{m-1}(x,\xi-i\epsilon\omega^{-1}\xim^{\kappa}\theta) $ is bounded from below by
\[
c\epsilon^2\omega^2\xim^2\prod_{k\neq i,j}|q_k(x,\xi-i\epsilon\omega^{-1}\xim^{\kappa}\theta)|^2.
\]
Summing up over all pair $i$, $j$ ($i\neq j$) we get the assertion for the case $j=2$. Continuing this argument one can prove   the case $j\geq 3$.
\qed

\medskip
Let us put
\[
h(x,\xi)=h_{m-1}(x,\xi-i\epsilon\omega^{-1}\xim^{\kappa}\theta)^{1/2}.
\]
\begin{lemma}
\label{lem:prebb}
There exists $C>0$ such that we have
\[
\left\{\begin{array}{ll}
|p^{(\al)}_{(\be)}|\leq C(\epsilon\omega)^{1-|\al+\be|}\xim^{1-|\al|}h,\quad 1\leq |\al+\be|\leq m,\\[3pt]
|p\,p^{(\al)}_{(\be)}|\leq C(\epsilon\omega)^{2-|\al+\be|}\xim^{2-|\al|}h^2,\quad 2\leq |\al+\be|\leq m.
\end{array}\right.
\]
\end{lemma}
\noindent
Proof: From \cite[Proposition 3]{Br} one has
\begin{align*}
|p^{(\al)}_{(\be)}(x,\xi)|\leq Ch_{m-|\al+\be|}(x,\xi)^{1/2}|\xi|^{|\be|}
\end{align*}
for $|\al+\be|\leq m$ which is bounded by $Ch_{m-|\al+\be|}(x,\xi-i\epsilon\omega^{-1}\xim^{\kappa}\theta)^{1/2}|\xi|^{|\be|}$ clearly.  On the other hand it follows from Lemma \ref{lem:preaa} 
\[
Ch(x,\xi)\geq (\epsilon\omega)^{|\al+\be|-1}\xim^{|\al+\be|-1}h_{m-|\al+\be|}(x,\xi-i\epsilon\omega^{-1}\xim^{\kappa}\theta)^{1/2}
\]
for $1\leq |\al+\be|\leq m$ which proves the assertion. The proof of the second assertion is similar.
\qed

\begin{lemma}
\label{lem:precc}Assume that $c^{\al}_{\be}\in S(\xim^{\rho|\be|-\delta|\al|},\bg)$ and $c^{\al}_{\be}\in S(\omega^{-1}\xim^{\kappa-|\al|},\bg)$ for $|\al+\be|=1$. Then 
for $1\leq |\al+\be|\leq m$ we have $p^{(\al)}_{(\be)}c^{\be}_{\al}\in  {\tilde S}(\omega^{-1}\xim^{\kappa}h,\bg)$ 
and $\epsilon^{|\al+\be|}|p^{(\al)}_{(\be)}c^{\be}_{\al}|\leq C\epsilon \omega^{-1}\xim^{\kappa}h$ with $C>0$ independent of $\epsilon$.
\end{lemma}
\noindent
Proof: Let $2\leq |\al+\be|\leq m$ then since $1=\kappa+2\delta$ we see by Lemma \ref{lem:prebb}
\begin{align*}
\epsilon^{|\al+\be|}|p^{(\al)}_{(\be)}c^{\be}_{\al}|\leq C\epsilon\omega^{1-|\al+\be|}\xim^{1-|\al|+\rho|\al|-\delta|\be|}h\\
\leq C\epsilon\omega^{-1}\xim^{\kappa}\big(\omega^{-1}\xim^{-\delta})^{|\al+\be|-2}h\leq C\epsilon\omega^{-1}\xim^{\kappa}h.
\end{align*}
When $|\al+\be|=1$ noting $c^{\be}_{\al}\in S(\omega^{-1}\xim^{\kappa-|\be|},\bg)$ we get the same assertion. We next estimate 
$\sum p^{(\al+\mu')}_{(\be+\nu')}(c^{\be}_{\al})^{(\mu'')}_{(\nu'')}$. If $|\al+\mu'+\be+\nu'|\geq m$ we have
\begin{align*}
\big|p^{(\al+\mu')}_{(\be+\nu')}\omega^{\be(\mu'')}_{\al(\nu'')}\big|\leq \xim^{m-|\al+\mu'|}\xim^{\rho|\al|-\delta|\be|}\xim^{-\rho|\mu''|+\delta|\nu''|}
\\
\leq C_{\epsilon}\omega^{-(m-1)}h\xim^{1-|\al+\mu'|}\xim^{\rho|\al|-\delta|\be|+\rho|\mu'|-\delta|\nu'|}\xim^{-\rho|\mu|+\delta|\nu|}\\
\leq C_{\epsilon}\omega^{-1}\xim^{\kappa}(\omega^{-(m-2)}\xim^{-\delta(|\al+\mu'+\be+\nu'|-2)})\xim^{-\rho|\mu|+\delta|\nu|}h\end{align*}
where the right-hand side is bounded by $C_{\epsilon}
\omega^{-1}\xim^{\kappa-\rho|\mu|+\delta|\nu|}h$. We turn to the case $|\al+\mu'+\be+\nu'|\leq m$. From Lemma \ref{lem:prebb} it follows 
\begin{align*}
|p^{(\al+\mu')}_{(\be+\nu')}|\leq C_{\epsilon}\omega^{1-|\al+\mu'+\be+\nu'|}\xim^{1-|\al+\mu'|}h\\
\leq C_{\epsilon}(\omega^{-1}\xim^{-\delta})^{|\mu'+\nu'|}\omega^{1-|\al+\be|}\xim^{1-|\al|}h\xim^{-\rho|\mu'|+\delta|\nu'|}\end{align*}
therefor for $|\al+\be|\geq 2$ we see easily
\begin{align*}
\big|p^{(\al+\mu')}_{(\be+\nu')}c^{\be(\mu'')}_{\al(\nu'')}\big|
\leq C_{\epsilon}\big(\omega^{-1}\xim^{-\delta}\big)^{|\al+\be|-2}\omega^{-1}\xim^{\kappa}h\xim^{-\rho|\mu|+\delta|\nu|}\\
\leq C_{\epsilon}\omega^{-1}\xim^{\kappa}h\xim^{-\rho|\mu|+\delta|\nu|}
\end{align*}
which also holds for $|\al+\be|=1$ because $c^{\be}_{\al}\in S(\omega^{-1}\xim^{\kappa-|\be|},\bg)$.
\qed

\subsection{Definition of $Q(z)$ which separates $P_{\psi}(z)$}
\label{Sec:pTM}

We follow the arguments in \cite{KNW}. Let us define ${\tilde p}(x+iy,\xi+i\eta)$ by
\[
{\tilde p}(x+iy,\xi+i\eta)=\sum_{|\al+\be|\leq m}\frac{1}{\al!\be!}\dif_x^{\al}\dif_{\xi}^{\be}p(x,\xi)(iy)^{\al}(i\eta)^{\be}.
\]
Then $p_{\psi}$ given by Theorem \ref{thm:matome:6} is expressed as ${\tilde p}(z-i\epsilon H_{\psi})$, which one can also write as
$$
{\tilde p}(z-i\epsilon H_{\psi})=\sum_{j=0}^{m}\big(i\frac{\dif}{\dif
t}\big)^jp(z-\epsilon tH_{\psi})/j!\big|_{t=0}.
$$
Using this expression we define $Q(z)$ which separates ${\tilde p}(z-i\epsilon H_{\psi})$ by
$$
Q(z)=\epsilon^{-1}|{\tilde H}_{\psi}|^{-1}\big(\frac{\dif}{\dif t}\big)\sum_{j=0}^{m}\big(i\frac{\dif}{\dif
t}\big)^jp(z-\epsilon tH_{\psi})/j!\big|_{t=0}
$$
where  $
{\tilde H}_{\psi}=(\lr{\xi}_{\gamma}\nabla_{\xi}\psi,-\nabla_x\psi)$. By the homogeneity it is clear that
\[
{\tilde p}(z-i\epsilon H_{\psi})
=\lr{\xi}_{\gamma}^m{\tilde p}({\tilde z}-i \lambda(z) 
{\tilde H}_{\psi}/|{\tilde H}_{\psi}|),\;\lambda(x,\xi)=\epsilon |{\tilde H}_{\psi}|\lr{\xi}_{\gamma}^{-1}
\]
where ${\tilde z}=(x,\xi\lr{\xi}_{\gamma}^{-1})$. It is not difficult to check  ${\tilde p}(z-i\epsilon H_{\psi})\in S(\lr{\xi}^m_{\gamma},\bg)$ and $Q\in S(\lr{\xi}^{m-1}_{\gamma},\bg)$. We study ${\tilde p}(z-i\epsilon H_{\psi})$ and $Q(z)$ in a conic neighborhood of $\rho$. We first recall 
\begin{pro}[{\cite[Lemma 5.8]{KW}}]
\label{pro:usimado}Let $\rho$ be a characteristic of $p$ of order $m$ and let $K\subset\Gamma_{\rho}$ be a compact set. Then one can find a conic neighborhood $V$ of $\rho$ and positive $C>0$ such that for any $(x,\xi)\in V$，$\zeta\in K$ and small $s\in\R$ one can write 
\[
p(z-s\zeta)=e(z,\zeta,s)\prod_{j=1}^m(s-\mu_j(z,\zeta))
\]
where $\mu_j(z,\zeta)$ are real valued and $e(z,\zeta,s)\neq 0$ for $(z,\zeta,s)\in V\times K\times \{|s|<s_0\}$. Moreover there exists $C>0$ such that we have 
\begin{equation}
\label{eq:sagawa}
|\mu_j(z,(0,\theta))|/C\leq |\mu_j(z,\zeta)|\leq C|\mu_j(z,(0,\theta))|,\;\;j=1,2,\ldots,m
\end{equation}
for any $(x,\xi)\in V$, $\zeta\in K$. 
\end{pro}
Writing $\prod_{j=1}^m(t-\mu_j)=\sum_{\ell=0}^mp_{\ell}\,t^{\ell}$ we see that ${\tilde p}(z-is\zeta)$ is written
\begin{eqnarray*}
\sum_{j=0}^{m}\frac{1}{j!}\big(is\frac{\dif}{\dif
t}\big)^j\big(e\prod_{j=1}^m(t-\mu_j)\big)\big|_{t=0}
=\sum_{\ell=0}^m\sum_{k=0}^{m-\ell}\frac{1}{k!}(is\frac{\dif}{\dif
t})^ke\big|_{t=0}p_{\ell}(is)^{\ell}
\end{eqnarray*}
which is equal to
\begin{eqnarray*}
\sum_{\ell=0}^m\bigl(\sum_{k=0}^{m}\frac{1}{k!}(is\frac{\dif}{\dif
t})^ke\big|_{t=0}-\sum_{
m-\ell+1\leq k\leq m}\frac{1}{k!}(is\frac{\dif}{\dif
t})^ke\big|_{t=0}\bigr)p_{\ell}(is)^{\ell}\\
=\sum_{k=0}^{m}\frac{1}{k!}(is\frac{\dif}{\dif
t})^ke\big|_{t=0}\prod_{j=1}^m(is-\mu_j)+O(s^{m+1})
\end{eqnarray*}
which proves
\begin{equation}
\label{eq:kehi}
{\tilde p}(z-is\zeta)
=e_{0}(z,\zeta,s)\prod_{k=1}^m(is-\mu_j(z,\zeta))+O(s^{m+1}).
\end{equation}
Note that $e_{0}(z,\zeta,s)=\sum_{k=0}^{m}(is\dif/\dif
t)^ke(z,\zeta,t)/k!\big|_{t=0}
$ and hence we have $e_{0}(z,\zeta,0)=e(z,\zeta,0)\neq 0$.

\begin{lemma}
\label{levithreei}There exist  a conic neighborhood $U$ of ${\rho}$ and a compact convex set $K\subset \Gamma_{{\rho}}$ such that ${\tilde H}_{\psi}/|{\tilde H}_{\psi}|\in K$ for $(x,\xi)\in U$, $\gamma\geq \gamma_0$.
\end{lemma}
\noindent
Proof: From $\psi=\xim^{\kappa}\log{(\phi+\omega)}$ it is easy to see
\[
\left\{\begin{array}{ll}\nabla_x\psi=\omega^{-1}\xim^{\kappa}\nabla_x\phi,\\
\nabla_{\xi}\psi=\omega^{-1}\xim^{\kappa}\nabla_{\xi}\phi+O(\xim^{\kappa-1})\log{(\phi+\omega)}+O(\xim^{\kappa-1}).\end{array}\right.
\]
Therefore one has
\begin{align*}
|{\tilde H}_{\psi}|^2
{\tilde H}_{\psi}=\omega^{-1}\xim^{\kappa}\big({\tilde H}_{\phi}+(O(1)\omega\log{(\phi+\omega)},0)\big).
\end{align*}
Since $|\phi+\omega|\leq 2\omega$ we can assume $\omega \log{(\phi+\omega)}$ is enough small taking $U$ small. In particular we have $\omega^{-1}\xim^{\kappa}/C\leq |{\tilde H}_{\psi}|\leq C\omega^{-1}\xim^{\kappa}$. Then noting  ${\tilde H}_{\phi}(\rho)/|{\tilde H}_{\phi}(\rho)|\in \Gamma_{\rho}$ which follows from  \eqref{eq:kuri} we get the assertion.
\qed

\medskip

We rewrite $Q(z)$ according to \eqref{eq:kehi}.
\begin{lemma}
\label{levifoura} 
Let 
${\tilde \omega}={\tilde H}_{\psi}/|{\tilde H}_{\psi}|$. Then we have
\begin{eqnarray*}
&Q(z)=\lr{\xi}_{\gamma}^{m-1}\bigl\{-i\dif e_{0}({\tilde z},{\tilde \omega},\lambda)/\dif \lambda
\prod_{j=1}^m(i\lambda-\mu_j({\tilde z},{\tilde \omega}))\\
&+e_{0}({\tilde z},{\tilde \omega},\lambda)\sum_{j=1
}^m\prod_{k=1,k\neq j}^m(i\lambda-\mu_{k}({\tilde z},{\tilde \omega}))+O(\lambda^m)\bigr\}.
\end{eqnarray*}
\end{lemma}
\noindent
Proof: Noting $\lambda(z)|{\tilde H_{\psi}}|^{-1}\lr{\xi}_{\gamma}=\epsilon$ one can write
\begin{align*}
Q(z)=
\epsilon^{-1}|{\tilde H}_{\psi}|^{-1}\lr{\xi}_{\gamma}^{m}\sum_{j=0}^{m}\big(\frac{\dif}{\dif
t}\big)\big(i\frac{\dif}{\dif t}\big)^jp({\tilde z}- t\lambda(z){\tilde \omega})/j!\big|_{t=0}\\
=\lr{\xi}_{\gamma}^{m-1}\sum_{j=0}^{m}\lambda(z)^j\big(\frac{\dif}{\dif
t}\big)\big(i\frac{\dif}{\dif t}\big)^jp({\tilde z}- t{\tilde \omega})/j!\big|_{t=0}
\end{align*}
which is equal to
\begin{equation}
\label{eqnvifoura}
\begin{split}
\lr{\xi}_{\gamma}^{m-1}\frac{\dif}{\dif s}\frac{1}{i}\sum_{j=0}^{m}
\frac{1}{(j+1)!}\big(is\frac{\dif}{\dif
t}\big)^{j+1}p({\tilde z}- t{\tilde \omega})\big|_{t=0,s=\lambda(z)}\\
=\frac{1}{i}\lr{\xi}_{\gamma}^{m-1}\frac{\dif}{\dif s}
\bigl\{{\tilde p}({\tilde z}-i s{\tilde \omega})-p({\tilde z})
+O(s^{m+1})\bigr\}_{s=\lambda(z)}\\
=\frac{1}{i}\lr{\xi}_{\gamma}^{m-1}
\bigl\{\frac{\dif}{\dif s}{\tilde p}({\tilde z}-i s{\tilde \omega})\big|_{s=\lambda(z)}
+O(\lambda^{m})\bigr\}.
\end{split}
\end{equation}
From \eqref{eq:kehi} the right-hand side of 
\eqref{eqnvifoura} turns to be
\begin{equation}
\label{eqnvifouri}
\begin{split}
\frac{1}{i}\lr{\xi}_{\gamma}^{m-1}\bigl\{\frac{\dif}{\dif s}\big(e_{0}({\tilde z},{\tilde \omega}, s)\prod_{j=1}
^m(i s-\mu_j({\tilde z},{\tilde \omega}))+O( s^{m+1})\big)\vert_{s=\lambda}
\bigr\}
\end{split}
\end{equation}
modulo $O(\lambda^{m})\xim^{m-1}$ which proves the assertion.
\qed
%

\subsection{Microlocal energy estimates}

To derive microlocal energy estimates we study ${\mathsf{Im}}\,(P_{\psi}\chi u,Q\chi u)$ where $\chi$ is a cutoff symbol supported in a conic neighborhood of $\rho$. Thus we are led to consider ${\mathsf{Im}}\,(P_{\psi}{\bar Q})$ in a conic neighborhood of $\rho$. Recall that $P_{\psi}=p_{\psi}+\sum_{j=0}^{m-1} (P_j)_{\psi}$ where $P_j(x,D)$ is the homogeneous part of degree $j$ of $P$ and $(P_j)_{\psi}\in S(\xim^j,\bg)$ by Theorem \ref{thm:matome:6}. Take any small $0<\epsilon^*\ll 1$ and we fix $\epsilon^*$ and put
\[
\delta=(1-\epsilon^*)/m,\;\; \rho=(m-1+\epsilon^*)/m,\;\;\kappa=\rho-\delta.
\]
\begin{lemma}
\label{lem:levifouri}Let $
S_0(z)={\mathsf{Im}} \big({\tilde p}(z-i\epsilon H_{\psi})\overline {Q(z)}\big)
$. Then one can find a conic neighborhood $V$ of ${\rho}$ and $C>0$ such that we have
\begin{align*}
\epsilon\, \omega^{-1}\xim^{\kappa}h^2(z)/C\leq S_0(z) 
\leq C\epsilon\,\omega^{-1}\xim^{\kappa}h^2(z).
\end{align*}
\end{lemma}
\noindent
Proof: Write $e_{0}({\tilde z},{\tilde \omega},\lambda)=e({\tilde z},{\tilde \omega},0)+i\lambda(\dif
e/\dif\lambda)({\tilde z},{\tilde \omega},0)+O(\lambda^2)$ then it is clear $
|e_{0}|^2=|e({\tilde z},{\tilde \omega},0)|^2+O(\lambda^2)$ so that $
{\mathsf{Re}}\,(\dif{\bar e}_{0}/\dif\lambda )e_{0}=2^{-1}\dif |e_{0}|^2/\dif\lambda
=O(\lambda)$. Thus from Lemma \ref{levifoura} and \eqref{eq:kehi} it follows that
\begin{eqnarray*}
&{\mathsf{Im}}\big({\tilde p}(z-i\epsilon H_{\psi})\overline{Q(z)}\big)
=\lr{\xi}_{\gamma}^{2m-1}|e_{0}|^2\lambda\\
&\times \bigl\{\sum_{j=1}^m\prod_{k=1,k\neq j}^m(\lambda^2+\mu_k^2)
\big(1+O(\lambda+\sum_{j=1}^m|\mu_j|)\big)\bigr\}.
\end{eqnarray*}
Since $\mu_j(\rho,{\tilde \omega})=0$, $j=1,2,\ldots,m$ one obtains
$$
S_0(z)\approx\lr{\xi}_{\gamma}^{2m-1}\lambda\sum_{j=1}^m\prod_{k=1,k\neq j}^m\big(\lambda^2+\mu_k({\tilde z},{\tilde \omega})^2\big).
$$
On the other hand noting $\lambda(z)\approx \epsilon\,\omega^{-1}\xim^{\kappa-1}$ and
\begin{align*}
&h_{m-1}(x,\xi-i\epsilon \omega^{-1}\xim^{\kappa}\theta)\\
&\approx\lr{\xi}_{\gamma}^{2m-2}\sum_{j=1}^m\prod_{k=1,k\neq j}^m\big(\epsilon^2\omega^{-2}\xim^{2\kappa-2}+\mu_k({\tilde z},(0,\theta))^2\big)
\end{align*}
we conclude the assertion from Lemma \ref{levithreei} and Proposition \ref{pro:usimado}.
\qed
\begin{lemma}
\label{leviithreea} We have $Q\in {\tilde S}(h,\bg)
$ and $S_0^{\pm 1}\in {\tilde S}((\omega^{-1}\xim^{\kappa}h^2)^{\pm 1},\bg)$ in $U$ for $\gamma\geq \gamma_0(\epsilon)$. Moreover $|Q|\leq C h$ with $C>0$ independent of $\epsilon$.
\end{lemma}
\noindent
Proof: From the definition one can see easily that $Q$ is a sum of terms, up to constant factor;
\begin{equation}
\label{eq:shaken}
\epsilon^{|\al+\be|-1}p^{(\al)}_{(\be)}(z)(\nabla_{\xi}\psi)^{\be}(\nabla_x\psi)^{\al}/(\lr{\xi}_{\gamma}^2|\nabla_{\xi}\psi|^2
+|\nabla_x\psi|^2)^{1/2} 
\end{equation}
with $1\leq|\al+\be|\leq m+1$. We also note that ${\mathsf{Im}}\,Q$ is a sum of such terms \eqref{eq:shaken} with $2\leq |\al+\be|\leq m+1$. From Lemma \ref{lem:logphi} it follows  
\[
\left\{\begin{array}{ll}
\nabla_{\xi}\psi\in S(\omega^{-1}\xim^{\kappa-1}, \bg),\\
\nabla_x\psi\in S(\omega^{-1}\xim^{\kappa}, \bg)
\end{array}\right.
\]
in $U$ and hence 
$(\lr{\xi}_{\gamma}^2|\nabla_{\xi}\psi|^2+|\nabla_x\psi|^2)^{-1/2}\in S((\omega^{-1}\xim^{\kappa})^{-1},\bg)$ then we have
\begin{align*}
V^{\be}_{\al}=(\nabla_{\xi}\psi)^{\be}(\nabla_x\psi)^{\al}(\lr{\xi}_{\gamma}^2|\nabla_{\xi}\psi|^2+|\nabla_x\psi|^2)^{-1/2}\\
\in
 S(\omega\xim^{-\kappa}\xim^{\rho|\al|-\delta|\be|},\bg).
\end{align*}
Noting $\omega^{-2}\xim^{\kappa-1}\leq 1$ it suffices to repeat the proof of Lemma \ref{lem:precc} to conclude $p^{(\al)}_{(\be)}V^{\be}_{\al}\in {\tilde S}(h,\bg)$ and $\epsilon^{|\al+\be|-1}|p^{(\al)}_{(\be)}V^{\be}_{\al}|\leq C h$ with $C$ independent of $\epsilon$ for $1\leq|\al+\be|\leq m+1$. From Lemma \ref{lem:precc} it follows that  ${\tilde p}(z-i\epsilon H_{\psi})-p(z)\in {\tilde S}(\omega^{-1}\xim^{\kappa}h,\bg)$.  Since from Lemma \ref{lem:prebb} one can  check  that $p(z){\mathsf{Im}}\,Q\in {\tilde S}(\omega^{-1}\xim^{\kappa}h^2,\bg)$ we get the assertion for $S_0$.  The assertion for $S_0^{-1}$ follows from Lemma \ref{lem:levifouri} and $S_0S_0^{-1}=1$.
\qed

\medskip
From Theorem \ref{thm:matome:6} and Lemma \ref{lem:precc} one can write
\[
p_{\psi}-{\tilde p}(z-i\epsilon H_{\psi})=\sqrt{\epsilon}\, r+r_0
\]
where $r\in {\tilde S}(\omega^{-1}\xim^{\kappa}h,\bg)$ with $|r|\leq C\epsilon \omega^{-1}\xim^{\kappa}h$ and $r_0\in S(\xim^{m-\delta(m+1)},\bg)$. Thus $r\#{\bar Q}\in {\tilde S}(\omega^{-1}\xim^{\kappa}h^2,\bg)$  and $|r{\bar Q}|\leq C\epsilon \omega^{-1}\xim^{\kappa}h^2$. On the other hand  Lemma \ref{lem:preaa} shows 
\[
\xim^{m-\delta(m+1)}=\xim^{m-1+\kappa-\delta(m-1)}\leq C\epsilon^{1-m}\gamma^{-\delta}\omega^{-1}\xim^{\kappa}h
\]
so that we see $r_0\in S(\omega^{-1}\xim^{\kappa}h,\bg)$ and $|r_0|\leq C\epsilon^{3/2}\omega^{-1}\xim^{\kappa}h$ for $\gamma\geq \gamma_0(\epsilon)$. Therefore in virtue of Lemma \ref{leviithreea} there is ${\tilde r}\in {\tilde S}(1,\bg)$ with $|{\tilde r}|\leq C\sqrt{\epsilon}$ such that 
\[
{\mathsf{Im}}\,(p_{\psi}\#{\bar Q})=S_0(1-{\tilde r})
\]
in some conic neighborhood of $\rho$.

We turn to $R=\sum_{j=0}^{m-1}(P_j)_{\psi}\in S(\xim^{m-1},\bg)$. 
From Lemma \ref{lem:preaa} again we have
\begin{align*}
\xim^{m-1}\leq C\epsilon^{1-m}
\omega^{-(m-1)}h
=C\epsilon^{1-m} \omega^{-1}\xim^{\kappa}h(\omega^{-(m-2)}\xim^{-\kappa})\\
\leq C\epsilon^{1-m}\gamma^{-\kappa+\delta(m-2)}\omega^{-1}\xim^{\kappa}h.
\end{align*}
Recalling that $\kappa-\delta(m-2)=\epsilon^*>0$ and hence we can assume  $C\epsilon^{1-m}\gamma^{-\epsilon^*}\leq C\epsilon^{3/2}$ for $\gamma\geq \gamma_0(\epsilon)$ so that there exists ${\hat r}\in {\tilde S}(1,\bg)$ with $|{\hat r}|\leq C\sqrt{\epsilon}$ such that $
{\mathsf{Im}}\,(R\#{\bar Q})=S_0(1-{\hat r})
$ in a conic neighborhood of $\rho$. Thus we conclude 
\[
\left\{\begin{array}{ll}
{\mathsf{Im}}\,(P_{\psi}\#{\bar Q})=E^2,\;\;E\in {\tilde S}(\omega^{-1/2}\xim^{\kappa/2}h,\bg),\\[3pt]
\epsilon^{1/2}\omega^{-1/2}\xim^{\kappa/2}h/C\leq |E|\leq C\epsilon^{1/2}\omega^{-1/2}\xim^{\kappa/2}h\\
\end{array}\right.
\]
in a conic neighborhood of $\rho$ with $C$ independent of $\epsilon$. The rest of the proof of deriving microlocal energy estimates is just a repetition of the arguments in \cite{KNW} and we conclude that the Cauchy problem for $p+P_{m-1}+\cdots$ is $\gamma^{(1/\kappa)}$ well-posed at the origin. Note that $1/\kappa=m/(m-2+2\epsilon^*)$ and $\epsilon^*>0$ is arbitrarily small so that $1/\kappa$ is as close to $m/(m-2)$ as we please.


\noindent
Department of Mathematics, Osaka University, 
Machikaneyama 1-1, Toyonaka, 560-0043, Osaka, Japan;
 nishitani@math.sci.osaka-u.ac.jp


\begin{thebibliography}{99}


\bibitem{Bea}
{R.Beals:} {\it Characterization of pseudodifferential operators and applications}, Duke Math. J. {\bf 44} (1977), 45--57.


\bibitem{BN4}
{ E.Bernardi and T.Nishitani:} {\it On the Cauchy problem for noneffectively hyperbolic operators: The Gevrey 4 well-posedness}, 
Kyoto J. Math. {\bf 51} (2011), 767--810.

\bibitem{BN3}
{ E.Bernardi and T.Nishitani:} {\it  On the Cauchy problem for noneffectively hyperbolic operators. The Gevrey 3 well-posedness}, 
J. Hyperbolic Differ. Equ. {\bf 8} (2011), 615--650.

\bibitem{Br}
{M.D.Bronshtein:} {\it The Cauchy problem for hyperbolic operators with characteristics of variable multiplicity}, Trudy  Moskov Mat. Obsc. {\bf 41} (1980), 83--99.

\bibitem{CT}
{F.Colombini and G.Taglialatela:} {\it Well-posedness for hyperbolic higher order operators with finite degeneracy}, J. Math. Kyoto Univ. {\bf 46} (2006), 833--877.

\bibitem{Ho2}
{ L.H\"ormander:} {\it The Cauchy problem for differential equations with double characteristics,} 
J. Anal. Math. {\bf 32} (1977), 118-196.

\bibitem{Hobook:iti}
{L.H\"ormander:} The analysis of linear 
partial differential operators. I. Distribution theory and Fourier analysis, Springer, 1990.

\bibitem{Hobook}
{ L.H\"ormander:} The analysis of linear partial differential operators. III. Pseudo-differential operators, Springer-Verlag, Berlin, 1994.

\bibitem{IP}
{ V.Ja.Ivrii and V.M.Petkov:} {\it Necessary conditions for the Cauchy problem for non strictly hyperbolic equations to be well posed,}  Uspehi Mat. Nauk. {\bf 29} (1974), 3-70.

\bibitem{Iv0}
{ V.Ja.Ivrii:} {\it Well-posedness conditions in Gevrey classes for the Cauchy problem for hyperbolic operators with characteristics of variable multiplicity,} Sib. Math. J. {\bf 17} (1977), 921-931.


\bibitem{Iwa2}
{ N.Iwasaki:} {\it The Cauchy problem for effectively hyperbolic equations (standard type),}
Publ. RIMS Kyoto Univ. {\bf 20} (1984), 551-592.


\bibitem{KN}
{ K.Kajitani and T.Nishitani:} The hyperbolic Cauchy problem, Lecture Notes in Mathematics, {\bf 1505}, Springer-Verlag, Berlin, 1991.

\bibitem{KNW}
{ K.Kajitani, S.Wakabayashi and T.Nishitani:} {\it The Cauchy problem for hyperbolic operators of strong type,} Duke Math. J. {\bf 75} (1994), 353-408.

\bibitem{KW}
{ K.Kajitani and S.Wakabayashi:} {\it Microlocal a priori estimates and the Cauchy problem. II,} Japan J. Math. {\bf 20} (1994), 1-71.

\bibitem{Ler}
{ N.Lerner:} Metrics on the phase space and non-selfadjoint pseudo-differential operators. Pseudo-Differential Operators. Theory and Applications, 3. Birkh\"auser Verlag, Basel, 2010.

\bibitem{Ni0}
{T.Nishitani:} {\it On the Lax-Mizohata theorem in the analytic and Gevrey classes,} J. Math. Kyoto Univ. {\bf 18} (1978), 509-521.

\bibitem{Ni11}
{ T.Nishitani:} {\it Local energy integrals for effectively hyperbolic operators I, II,} 
J. Math. Kyoto Univ. {\bf 24} (1984), 623-658 and 659-666.

\bibitem{Ni:2}
{T.Nishitani:} {\it Necessary conditions for strong hyperbolicity of first order systems,} 
J. Anal. Math. {\bf 61} (1993), 181-229.


\bibitem{NTa}
{ T.Nishitani and M.Tamura:} {\it A class of Fourier integral operators with complex phase related to the Gevrey classes,} J. Pseudo-Differ. Oper. Appl. {\bf 1} (2010), 255-292.


\bibitem{Tez}
{A. M. Tezin:} {\it Integral curves of a generalized homogeneous differential equations of the first order,} Differential Equations {\bf 1} (1965), 570-580.




\end{thebibliography}
\end{document}